\newif\ifarxiv
\arxivtrue 

\documentclass[ijoc,nonblindrev]{informs3} 

\OneAndAHalfSpacedXII 
\usepackage{hyperref}

\usepackage{anyfontsize}
\usepackage{newtxtext,newtxmath,amsfonts,amssymb}
\usepackage{mathtools}
\usepackage{multirow}
\usepackage{nomencl}
\renewcommand{\nomgroup}[1]{%
\ifthenelse{\equal{#1}{S}}{\item[\textbf{Indices}]}{%
\ifthenelse{\equal{#1}{P}}{\item[\textbf{Parameters}]}{%
\ifthenelse{\equal{#1}{V}}{\item[\textbf{Variables}]}{%
\ifthenelse{\equal{#1}{R}}{\item[\textbf{Random Variables}]}{}}}}
}

\makenomenclature
\usepackage{makecell}
\usepackage{graphicx}
\PassOptionsToPackage{hyphens}{url}\usepackage{hyperref}
\usepackage[utf8]{inputenc}
\usepackage[english]{babel}
\usepackage{amsmath}
\usepackage{gensymb}
\usepackage{bbm}
\usepackage{amssymb}

\usepackage{tikz}
\tikzset{
    block/.style = {rectangle, draw, text width=4cm, text centered, minimum height=1.5cm, node distance=1.5cm},
    smallblock/.style = {rectangle, draw, text width=3.2cm, text centered, minimum height=1cm, node distance=1.5cm},
    line/.style = {draw, -latex'},
}
\usepackage{pgfplots}
\usetikzlibrary{intersections}
\usepgfplotslibrary{fillbetween}
\usetikzlibrary{arrows, arrows.meta, positioning, backgrounds, automata, decorations, shapes, fit, calc}

\DeclareFontFamily{U}{mathx}{\hyphenchar\font45}
\DeclareFontShape{U}{mathx}{m}{n}{
      <5> <6> <7> <8> <9> <10>
      <10.95> <12> <14.4> <17.28> <20.74> <24.88>
      mathx10
      }{}
\DeclareSymbolFont{mathx}{U}{mathx}{m}{n}
\DeclareFontSubstitution{U}{mathx}{m}{n}
\DeclareMathAccent{\widecheck}{0}{mathx}{"71}
\DeclareMathAccent{\wideparen}{0}{mathx}{"75}

\usepackage{amsfonts}
\usepackage{chngcntr}
\usepackage{apptools}
\usepackage{pgfgantt}
\usepackage{lscape}
\usepackage[pass]{geometry}
\usepackage{booktabs}
\usepackage{algpseudocode} 
\usepackage[caption=false, font=footnotesize, subrefformat=parens,labelformat=parens]{subfig}

\usepackage{algorithm}
\usepackage{array}

\usepackage{textcomp, gensymb}

\usepackage{xcolor}
\usepackage{cases} 
\usepackage{empheq}
\usepackage{enumitem}

\def\b#1{\boldsymbol{#1}}

\mathchardef\mhyphen="2D 
\usepackage{natbib}
 \bibpunct[, ]{(}{)}{,}{a}{}{,}%
 \def\bibfont{\small}%
\TheoremsNumberedThrough     

\EquationsNumberedThrough    

\pgfplotsset{compat=1.18}
\begin{document}


\setlength{\floatsep}{8pt plus 1pt minus 2pt}       



\TITLE{Strengthened and Faster Linear Approximation to Joint Chance Constraints with Wasserstein Ambiguity}

\ARTICLEAUTHORS{%
\AUTHOR{Yihong Zhou}
\AFF{Department of Engineering Science, University of Oxford, U.K., \EMAIL{yihong.zhou@eng.ox.ac.uk}, \URL{}}
\AUTHOR{Yuxin Xia}
\AFF{Department of Engineering Science, University of Oxford, U.K., \EMAIL{yuxin.xia@eng.ox.ac.uk}, \URL{}}
\AUTHOR{Hanbin Yang}
\AFF{Business School, Hong Kong University of Science and Technology, Hong Kong, China, \EMAIL{hanbinyang@ust.hk}, \URL{}}
\AUTHOR{Thomas Morstyn}
\AFF{Department of Engineering Science, University of Oxford, U.K., \EMAIL{thomas.morstyn@eng.ox.ac.uk}, \URL{}}
} 

\ABSTRACT{Many real-world decision-making problems in energy systems, transportation, and finance have uncertain parameters in constraints. Wasserstein distributionally robust joint chance constraints (WDRJCC) offer a promising solution by explicitly guaranteeing the probability of the simultaneous satisfaction of multiple constraints. However, WDRJCC are computationally demanding, and practical applications often require more tractable approaches, especially for large-scale and complex problems such as power system unit commitment problems and multilevel problems with chance constraints in lower levels. To address this, this paper proposes a novel convex inner-approximation for WDRJCC with right-hand-side uncertainties (RHS-WDRJCC). Motivated by the strengthening process that leads to a faster but still exact mixed-integer reformulation, we propose a Strengthened and Faster Linear Approximation (SFLA) by strengthening an existing convex inner-approximation. This strengthening process reduces the number of constraints and tightens the feasible region for ancillary variables, leading to significant computational speedup. We prove that the proposed SFLA does not introduce additional conservativeness and can even be less conservative compared to common approximations such as W-CVaR. \textcolor{black}{We then extend the proposed SFLA to robustness maximization, a decision-making paradigm that can be more interpretable, where the risk level and the Wasserstein radius are determined by maximizing solution robustness subject to a utility degradation limit. We discuss the connection between risk minimization and radius maximization as two formulations of robustness maximization, and show the advantage of radius maximization.}

In power system unit commitment, the proposed SFLA achieves up to $10\times$ and on average $3.8\times$ computational speedup compared to the strengthened and exact mixed-integer reformulation in finding comparable high-quality solutions. In a bilevel strategic bidding problem where the exact reformulation is not applicable due to non-convexity, the proposed SFLA can lead to $90\times$ speedup compared to existing convex approximation methods including W-CVaR. In robustness maximization, the proposed SFLA demonstrated over $100\times$ speedup than other convex-approximations.}
%



\maketitle

%


\section{Introduction}\label{sec:intro}


Many real-world decision-making problems in energy systems, transportation, and finance involve constraints with uncertain parameters due to variable supply, customer demand, and economic conditions \citep{gabrel2014recent}. A typical solution approach is robust optimization (RO), which aims to ensure feasibility for each possible scenario; however, this often leads to overly conservative solutions. In contrast, chance-constrained programming (CCP) provides a less restrictive alternative~\citep{Stochastic-programming}. CCP supports robust but economic decision-making by explicitly limiting the probability of constraint violations, allowing for a controlled degree of conservativeness relative to RO. As a result, CCP has found widespread applications in fields such as power systems, economics, finance, and water management~\citep{CCP_power_system_review}.

However, classic CCP models are based on the exact distribution of random variables, which is typically not available, and decision-makers generally only possess a historical dataset. This dataset may be insufficient to accurately infer the true distribution of random variables, which limits the out-of-sample \textcolor{black}{performance} of a CCP model. To hedge the ambiguity of the underlying distribution, distributionally robust chance-constrained programming (DRCCP) has been proposed to control the violation probability under the worst-case probability distribution over a so-called ``ambiguity'' set~\citep{scarf1957min}. This ambiguity set may be defined as the set of probability distributions with the same \textcolor{black}{statistical} moments, such as mean and variance~\citep{delage2010distributionally}, or defined as distributions within a certain distance from the reference distribution. Despite good tractability, moment-based ambiguity sets do not fully use the information of the datasets and may lead to over-conservativeness~\citep{gao2023distributionally}. For distance-based ambiguity sets, the Wasserstein distance is widely applied due to its superior out-of-sample performance~\citep{mohajerin2018data}, and the advantage of Wasserstein-based ambiguity sets over other distance-based sets was discussed by~\cite{gao2023distributionally}. A Wasserstein DRCCP model can be written as:
\begin{subequations}
\label{eq:wdrccp}
\begin{align}
    \min_{\b x \in \mathcal{X}} \quad & c(\b x) \label{obj:wdrjcc}  \\
    \mbox{s.t.} \quad & \sup_{\mathbb {P} \in \mathcal{F}_N (\theta)} \mathbb{P}[\b \xi \notin \mathcal{S}(\b x)] \leq \epsilon , \label{constr:wdrjcc}
\end{align}
\end{subequations}
where $c(\cdot)$ is the objective function, $\mathcal{X} \subset \mathbb{R}^L$ is a compact domain for the decision variables $\b x \in \mathbb{R}^L$. Constraint~\eqref{constr:wdrjcc} ensures that the random vector $\b \xi \in \mathbb{R}^K$ falls outside the decision-dependent safety set $\mathcal{S}(\b x)$ with a small probability, which is not greater than the risk level $\epsilon$ under the worst-case distribution $\mathbb{P}\in \mathcal{F}_N (\theta)$. Here, the Wasserstein ambiguity set $\mathcal{F}_N(\theta)$ is defined as a ball with radius $\theta$ centered at the empirical distribution $\mathbb{P}_N$ constructed by historical data $\mathcal{F}_N(\theta) \coloneqq \left\{ \mathbb{P} \mid  d_W(\mathbb{P}_N, \mathbb{P}) \leq \theta \right\}$.

Following \cite{2022_OR_exact_DRO} and \cite{exact_milp_strengthened}, we consider the $1$-Wasserstein distance based on a general norm $\|\cdot\|$, expressed as $d_W (\mathbb{P}, \mathbb{P}') \coloneqq \inf\limits_{\Pi \in \mathcal{P}(\mathbb{P}, \mathbb{P}')} \mathbb{E}_{(\b \xi, \b\xi') \sim \Pi} [\| \b\xi - \b\xi' \|]$, where $\mathcal{P}(\mathbb{P}, \mathbb{P}')$ is a set of distributions with marginal distributions $\mathbb{P}$ and $\mathbb{P}'$.
This paper focuses on Wasserstein distributionally robust joint chance constraints under right-hand-side uncertainty (RHS-WDRJCC)~\eqref{constr:wdrjcc} with the following safety set:
\begin{align}
\label{eq:def:safety_set}
    \mathcal{S}(\b x) \coloneqq \left\{ \b \xi \mid \b a_p^\top \b x \leq \b b_p^\top \b \xi + d_p, \  p\in [P] \right\},
\end{align}
where there are $P$ constraints indexed by $[P] \coloneqq \{1, \cdots, P\}$ that need to be met jointly with high probability $1-\epsilon$ with $\epsilon \in (0,1)$. This joint satisfaction is desired in practical applications for its higher safety~\citep{CCP_power_system_review, ding2022distributionally}, and the RHS uncertainty also arises in several practical problems, such as managing the thermal constraints of the power grid lines, securing the power system reserve~\citep{wang2016risk, wu2016solution, yang2019analytical}, or restricting the load-generation imbalance~\citep{van2018exact}. 

To exactly solve the problem with RHS-WDRJCC~\eqref{eq:wdrccp}, \cite{2022_OR_exact_DRO} proposed an exact mixed-integer programming (MIP) reformulation, which was further strengthened by~\cite{exact_milp_strengthened} by exploiting valid inequalities and reducing the values of big-M. 
\cite{jiang2024terminator} further proposed a method that combines inner and outer approximations to derive optimality cuts, thereby accelerating the exact MIP solution process. Alternative methods include ALSO-X, ALSO-X+, and ALSO-X\# that (approximately) solve a CCP in an iterative way with high quality~\citep{jiang2022also, jiang2024also}.
However, the MIP-based exact reformulation is NP-hard, and these ALSO methods may have numerical issues and \textcolor{black}{a} lack of stability due to their iterative schemes. Greater tractability is still desired, especially when the original deterministic problem is already complicated. We have identified two types of practical applications that necessitate more computationally efficient solution schemes.

One typical example is the power system unit commitment (UC) problem that determines the on/off statuses of generators to minimize operational cost within power network constraints, which is a large-scale multi-period MIP problem with thousands of generators and network nodes included. In one of the largest electricity markets in the world, managed by the Midcontinent Independent
System Operator (MISO), the network model includes more than $45{,}000$ buses, $1{,}400$ generation sources, around $2{,}500$ pricing nodes, and has a $36$-hour look-ahead horizon with hourly resolution \citep{MISO_UC_difficulty, sun2018novel}. Additionally, UC problems may be solved repeatedly during the market-clearing process to accommodate necessary modifications to constraints or other evaluations. Moreover, the entire process must be completed within a stringent time frame of three to four hours~\citep{MISO_UC_difficulty}. The UC complexity is further exacerbated by the increasing number of virtual bids (bids for speculation submitted by purely financial players with no obligation to own physical generation) and the increasing number of combined cycle units that require complicated operation models \citep{sun2018novel}. However, the increasing penetration of fluctuating renewables increases the uncertainty level of the system, and thus it becomes more important to have chance-constrained UC to ensure reliable operation ~\citep{yang2019analytical,van2018exact}.

The second example involves applying RHS-WDRJCC at lower levels of bilevel or multilevel optimization problems, where all levels need to reach the optimum. To solve these multilevel problems, Karush-Kuhn-Tucker (KKT) conditions or strong duality need to be exploited to derive solvable single-level counterparts, which require convexity in lower levels. Practical applications include strategic price-maker bidders in day-ahead energy markets \citep{paredes2023stacking} and gas markets that involve a sequential clearing process~\citep{heitsch2022convex}. In these applications, market clearing is modeled at lower levels. Because the clearing process can occur before uncertainties are revealed (e.g., day-ahead markets), a CCP needs to be implemented at these lower levels to manage uncertainties~\citep{beck2023survey}. However, the MIP-based RHS-WDRJCC reformulation \citep{exact_milp_strengthened} is unsuitable due to its non-convexity, and the iterative ALSO methods~\citep{jiang2022also, jiang2024also} are less capable of incorporating optimality conditions.

Due to the high complexity of these practical problems, it is desirable to have convex or linear approximations of RHS-WDRJCC, so as to bring a minimum extra computational effort while maintaining solution quality. Common approaches include the Bonferroni approximation~\citep{chen2023approximations}, which can become overly conservative when events overlap significantly, and the worst-case conditional value at risk (W-CVaR)~\citep{mohajerin2018data, chen2023approximations}, a common approximation for CCPs. \cite{chen2023approximations} also introduced the ``best'' convex inner-approximation for RHS-WDRJCC, which is equivalent to W-CVaR under certain hyperparameters. However, these approximation schemes introduce a large number of additional ancillary variables and constraints, which brings excessive computing burden, especially for the UC and multilevel problems discussed above; consequently, CCPs generally remain unscalable for industrial-scale instances~\citep{zhao2024viabilitystochasticeconomicdispatch}. Therefore, it is still preferable to further \textcolor{black}{enhance the} computational efficiency of RHS-WDRJCC.

\begin{figure}
\vspace{-4 mm}
     \FIGURE
    {
    \includegraphics[width=1\linewidth]{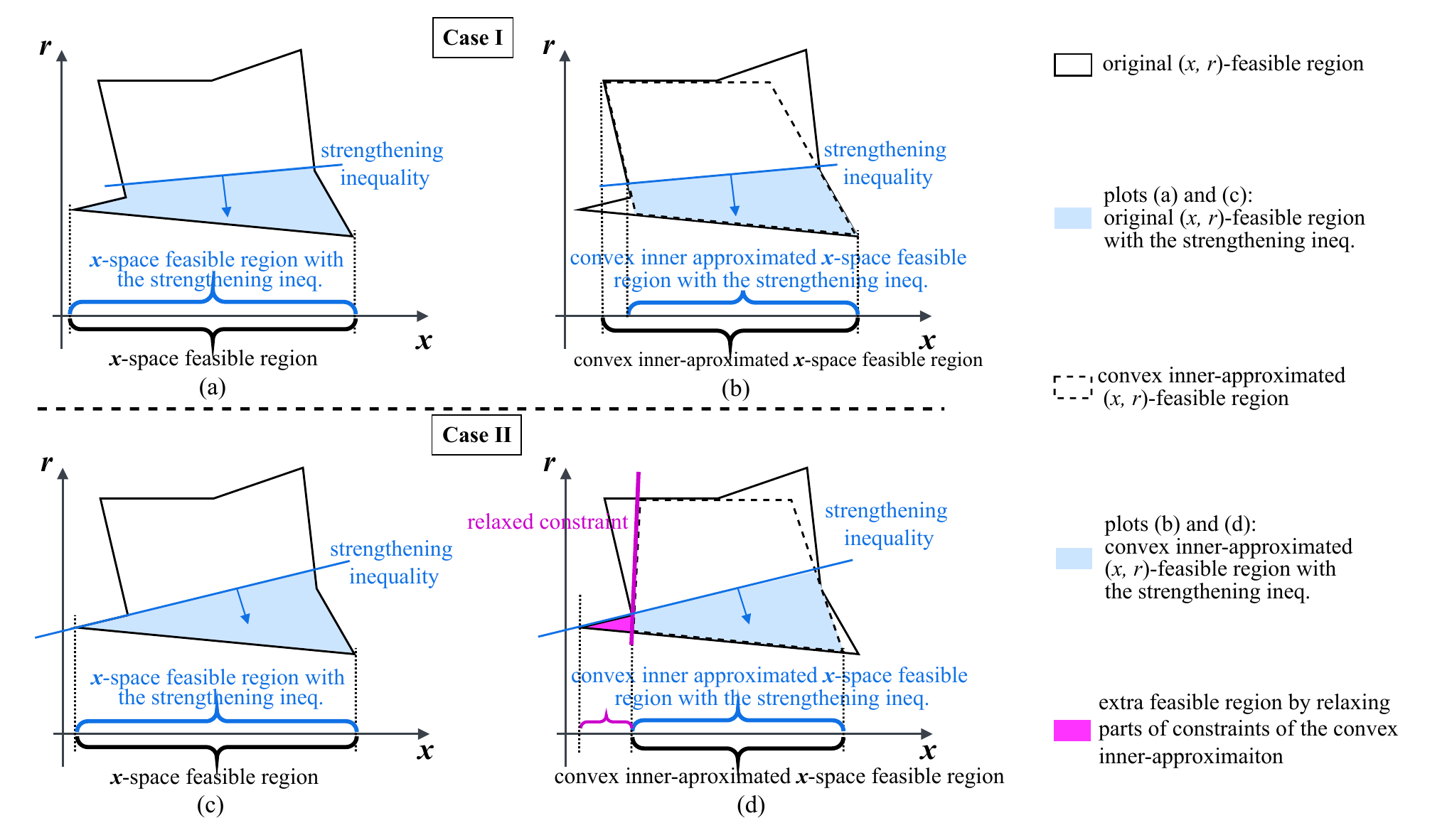}
    \vspace{-2 mm}
    } 
{Two possible cases when a strengthening inequality valid for the original feasible region is applied to the convex inner-approximation. \label{fig:valid_ineq}}
{Analytical reformulations sometimes (such as for WDRJCC) introduce \textcolor{black}{an} ancillary variable $r$ in addition to the original decision variable $x$. Case I: A strengthening inequality is valid for the original nonconvex region but leads to higher conservativeness \textcolor{black}{in its inner approximation}. (a) In the original nonconvex feasible region of $x$ and $r$, the strengthening (valid) inequality excludes a redundant portion without reducing the feasible region of $x$. (b) The convex inner-approximation (dashed) has a smaller feasible region of $x$ compared to the original region, and the strengthening inequality further reduces the $x$ feasible region, resulting in even higher conservativeness. Case II: A strengthening inequality is valid for both the original nonconvex region and its inner approximation. (c) In the original nonconvex feasible region, the strengthening inequality excludes a redundant portion without reducing the feasible region of $x$. (d) The convex inner-approximation (dashed) has a smaller feasible region of $x$ compared to the original region, and the strengthening inequality does not further reduce the $x$ feasible region. Furthermore, by relaxing certain constraints of the convex inner-approximation, it is possible to include an extra non-empty feasible region (depicted in purple) and thus expand the feasible region of $x$ (less conservativeness). This extra feasible region is small relative to the excluded region; therefore, the overall search space is still reduced, preserving the computational efficiency provided by the strengthening inequality.}
\end{figure}

One promising solution approach is to utilize valid inequalities, which are additional constraints introduced to an optimization problem to tighten the formulation without excluding any feasible solutions from the original problem's feasible region. Motivated by the work of \cite{exact_milp_strengthened}, where valid inequalities were exploited to strengthen the exact MIP reformulation so as to speed up computation, this paper shows that the strengthening process can also be migrated to a convex inner-approximation that admits W-CVaR equivalence. This migration reduces the number of constraints in the original convex approximation and \textcolor{black}{tightens} the feasible region for ancillary variables, \textcolor{black}{resulting in} significant computational speedups. \textcolor{black}{Although the migration may appear direct, there is a subtle source of conservativeness: migrating a strengthening process can introduce additional conservativeness when applied to a convex inner approximation, as shown in Case I of Figure~\ref{fig:valid_ineq}. However, and fortunately, we prove that the proposed SFLA will always lead to Case II in Figure~\ref{fig:valid_ineq}, which does not introduce extra conservativeness and can even produce a less conservative inner-approximation. The computational speedup and reduced conservativeness of the proposed SFLA address the practical needs of industry-scale problems, such as UC and bilevel strategic bidding.} In summary, the main contributions of this paper include:
\begin{enumerate}[label={(\arabic*)}, leftmargin=*, ref=(\arabic*), align=left]
    \item We proposed a Strengthened and Faster Linear Approximation (SFLA) to RHS-WDRJCC. By exploiting valid inequalities for a convex inner-approximation that is equivalent to W-CVaR~\citep{chen2023approximations} under specific hyperparameters, the number of constraints is reduced, and the feasible region of ancillary variables \textcolor{black}{is} tightened, which leads to significant computational speedup. Even with the tightening, we demonstrate that the proposed SFLA does not introduce additional conservativeness and can even provide a less conservative inner approximation. In addition, we offer theoretical guidance for setting the SFLA hyperparameter and analyze when SFLA becomes equivalent to exact reformulation or less conservative (or equivalent) to existing convex approximations such as W-CVaR. 

    \item We extend the proposed SFLA to robustness maximization, a decision-making paradigm that can be more interpretable to set the risk level and the Wasserstein radius. We show that the proposed SFLA applies to both risk minimization and radius maximization (the two formulations for robustness maximization), offering computational speedups without introducing extra conservativeness. We further show the equivalence of the feasible regions between risk minimization and radius maximization, and demonstrate that radius maximization can be both more computationally efficient and more robust.

    \item We demonstrated the superiority of the proposed SFLA through extensive numerical studies on two important real-world large-scale chance-constrained optimization problems. First, in the UC problem, the proposed SFLA achieves up to $10\times$ and on average $3.8\times$ computational speedup compared to the existing strengthened and exact MIP reformulation~\citep{exact_milp_strengthened} in finding comparable high-quality solutions. We also demonstrate that the proposed SFLA achieves approximation quality comparable to the exact reformulation, except for minor inferiority only occurring under the combination of a high risk level, a large number of historical data, and a small ambiguity set radius. Second, in the bilevel strategic bidding problem where the exact formulation is not applicable due to non-convexity, the proposed SFLA demonstrates $90\times$ speedup in computation time compared to convex inner-approximations such as W-CVaR. In robustness maximization, the proposed SFLA demonstrated over $100\times$ computational speedup compared to other convex approximations with the same solution quality.
\end{enumerate}
This paper is organized as follows. Section \ref{sec:problem_form} lists the existing exact reformulation and convex inner-approximations. Section \ref{sec:main_theo_res} introduces the proposed SFLA and demonstrates its theoretical properties that make it particularly useful for real-world large-scale problems. \textcolor{black}{Section~\ref{sec:robustsmax} extends SFLA to a more interpretable robustness-maximization framework.} Section \ref{sec:casestudies} carries out numerical studies and Section \ref{sec:conclusion} concludes this paper. 

\section{Problem Formulation}
\label{sec:problem_form}

This paper focuses on RHS-WDRJCC~\eqref{constr:wdrjcc}, with the safety set described in Eq.~\eqref{eq:def:safety_set}. In what follows, we assume $\epsilon \in (0,1)$ and $\theta > 0$, which are conditions necessary for the subsequent exact reformulation on which this work is based~\citep{2022_OR_exact_DRO}. This section summarizes the current exact reformulation and approximation methods for RHS-WDRJCC.

\subsection{Exact Reformulation}

Suppose we have collected $N$ independent and identically distributed (i.i.d.) samples $\{ \b \xi_i \}_{i\in [N]}$ for the random vector $\b \xi$, where the index set is defined as $[N] \coloneqq \{1, \cdots, N\}$. Given a decision $\b x \in \mathcal{X}$ and a sample $\b \xi_i$, we define the distance from $\b \xi_i$ to the complement of $\mathcal{S}(\b x)$ as:
\begin{equation}
    \mathrm{dist}\left(\b \xi_i, \mathcal{S}(\b x)\right) = \inf_{\b \xi' \in \mathbb{R}^K} \left\{ \| \b \xi_i - \b\xi' \| \mid  \b\xi' \notin \mathcal{S}(\b x) \right\}.
\end{equation}
Building on the results given by~\cite{2022_OR_exact_DRO} and the safety set defined in~\eqref{eq:def:safety_set}, we can reformulate $ \mathrm{dist}\left(\b \xi_i, \mathcal{S}(\b x)\right)$ to the following analytical expression:
\begin{equation}\label{eq:dist_anal}
    \mathrm{dist}\left(\b \xi_i, \mathcal{S}(\b x)\right)
        = \left(\min_{p\in [P]} \dfrac{\b b_p^\top \b \xi_i + d_p - \b a_p^\top \b x}{\| \b b_p \|_*} \right)^+ = \min_{p\in [P]} \left(\dfrac{\b b_p^\top \b \xi_i + d_p - \b a_p^\top \b x}{\| \b b_p \|_*} \right)^+,
\end{equation}
where $\|\cdot\|_*$ is the dual norm and $(\cdot)^+$ extracts the non-negative component of the argument, setting negative values to zero. By introducing ancillary variables $s\in \mathbb{R}$ and $\b r \in \mathbb{R}^N$, \cite{2022_OR_exact_DRO} showed that RHS-WDRJCC~\eqref{constr:wdrjcc} can be expressed as the following distance-based formulation:
\begin{subequations}\label{eq:wdrjcc_abstract}
    \begin{align}
        &  s\geq 0, \b r \geq \b 0, \label{eq:s_r} \\
        & \epsilon N s - \sum\limits_{i=1}^N r_i \geq \theta N, \label{eq:epsilonNs} \\
        & \left( \frac{\b b_p^\top \b \xi_i + d_p - \b a_p^\top \b x}{\| \b b_p \|_*} \right)^+ \geq s-r_i, &&\hspace{-50pt} \forall \ i\in [N], p \in [P]. \label{eq:dist_constr}
    \end{align}
\end{subequations}
While this distance-based formulation cannot be directly processed by commercial solvers such as Gurobi~\citep{gurobi}, it can be transformed into an equivalent MIP formulation:

This MIP reformulation, based on~\cite{2022_OR_exact_DRO}, is originally designed for open safety sets, but~\cite{exact_milp_strengthened} demonstrated that it applies to both open and closed sets. 

Despite being solvable by advanced MIP solvers, the MIP reformulation remains NP-hard~\citep{luedtke2010integer} and its complexity becomes particularly problematic when the original deterministic problem is already complicated or requires convexity for tractable reformulation such as in multilevel problems. Therefore, we focus on convex inner-approximation (for safety) of the $\b x$-feasible region $\mathcal{X}_\text{Exact} \coloneqq \left\{\b x \in \mathcal{X} \mid \exists s, \b r :~\eqref{eq:s_r}\mbox{--}\eqref{eq:dist_constr}\right\}$ defined by the RHS-WDRJCC in model~\eqref{eq:wdrccp} for the decision variable $\b x$,
because this is the region that affects the objective value $c(\b x)$ in~\eqref{obj:wdrjcc}.

\subsection{Convex Inner-Approximations}

In this section, we introduce three representative convex inner-approximation methods for RHS-WDRJCC.
The first is the Bonferroni approximation, which is a straightforward and computationally efficient \textcolor{black}{inner approximation} that replaces a joint chance constraint (JCC) with individual chance constraints at reduced individual risk levels $\epsilon_p$. According to~\cite{chen2023approximations}, the Bonferroni approximation of RHS-WDRJCC~\eqref{constr:wdrjcc} can be formulated as:
\begin{align}\label{eq:Bonfe_RHS}
    \b a_p^\top \b x \leq d_p - \sup_{\mathbb{P} \in \mathcal{F}(\theta)} \mathbb{P}\mhyphen\mathrm{VaR}_{\epsilon_p} (-\b b_p^\top \b \xi),\quad \forall p\in [P],
\end{align}
where the preset individual risk levels $\epsilon_p$ must satisfy $\sum_{p \in P} \epsilon_p \leq \epsilon$. Here, $\sup_{\mathbb{P} \in \mathcal{F}(\theta)} \mathbb{P}\mhyphen\mathrm{VaR}_{\epsilon_p} (-\b b_p^\top \b \xi)$ is the optimal value of a bilinear optimization problem~\citep{chen2023approximations} for a specific value of $\epsilon_p$.
The Bonferroni approximation~\eqref{eq:Bonfe_RHS} does not add extra variables or constraints to the original deterministic problem, resulting in minimal computational overhead. A common practice is to set $\epsilon_p = \epsilon / P$~\citep{CCP_power_system_review}, but this often leads to excessive conservativeness, especially when $P$ is large or when the safety constraints are correlated~\citep{chen2007robust}. However, even with optimized $\epsilon_p$, the approach can still be overly conservative in some cases, and finding optimal values for $\epsilon_p$ is generally intractable~\citep{chen2023approximations}. As will be shown in our case studies, the Bonferroni approximation often leads to infeasibility due to this conservativeness. Consequently, we exclude this approach from theoretical comparisons in the following discussion.

An alternative and widely applied inner-approximation is W-CVaR~\citep{mohajerin2018data, chen2023approximations}, which can be formulated as the following set of linear constraints:
\begin{subequations}\label{eq:CVaR_RHS}
    \begin{align}
        & \b \alpha \geq \b 0, \beta \in \mathbb{R}, \tau \in \mathbb{R}, \label{eq:wcvar-fea} \\
        & \tau + \frac{1}{\epsilon} \left(\theta \beta + \frac{1}{N} \sum_{i\in [N]} \alpha_i\right) \leq 0, \label{eq:wcvar-b}\\
        & \alpha_i \geq w_p \left(\b a_p^\top \b x - \b b_p^\top \b \xi_i - d_p\right) - \tau, &&\hspace{-50pt} \forall i\in[N], p\in [P],\label{eq:wcvar-c} \\
        & \beta \geq w_p \| \b b_p \|_*,  &&\hspace{-50pt} \forall p \in [P], \label{eq:wcvar-last}
    \end{align}
\end{subequations}
where $\b w \coloneqq [w_p]_{p\in[P]}$ subject to $\b w \in \Delta_{++}\coloneqq\{\b w \in (0,1)^P \mid \sum_{p\in[P]}w_p = 1\}$ is a tunable hyperparameter that affects the performance of the W-CVaR approximation by prioritizing specific constraints in $\mathcal{S}(\b x)$~\citep{chen2023approximations, ordoudis2021energy}. Similarly, we define its $\b x$-feasible region as $\mathcal{X}_\text{WCVaR}(\b w) \coloneqq
\{ \b x \in \mathcal{X} \mid \exists \b\alpha, \beta, \tau :~\eqref{eq:wcvar-fea}\mbox{--}\eqref{eq:wcvar-last} \}$.

Another convex approximation~\citep{chen2023approximations} \textcolor{black}{linearizes constraint~\eqref{eq:dist_constr}} by overcoming the non-convexity caused by the distance function~\eqref{eq:dist_anal}. \textcolor{black}{Specifically, it} replaces the original distance function $\text{dist}(\b \xi_i, \mathcal{S}(\b x))$ with a conservative approximation $\widehat{\text{dist}}(\b \xi_i, \mathcal{S}(\b x))$ given by:
\begin{align}\label{eq:dist_approx}
    \widehat{\text{dist}}(\b \xi_i, \mathcal{S}(\b x)) \coloneqq \kappa_i \left(\min_{p \in [P]} \frac{\b b_p^\top \b \xi_i + d_p - \b a_p^\top \b x}{\| \b b_p \|_*} \right) ,
\end{align}
where $\kappa_i \in [0, 1]$ are pre-selected parameters. Since $\widehat{\text{dist}}(\b \xi_i, \mathcal{S}(\b x)) \leq \text{dist}(\b \xi_i, \mathcal{S}(\b x))$, the replacement yields the following linear inner-approximation (LA) of the constraints in~\eqref{eq:wdrjcc_abstract}:
\begin{subequations}\label{eq:LP_constr_1}
    \begin{align}
        & s\geq 0, \b r \geq \b 0, \label{eq:LP1_s_r} \\
        & \epsilon N s - \sum\limits_{i\in [N]}r_i \geq \theta N,  \label{eq:LP1_epsilonNs}\\
        & \kappa_i \left( \frac{\b b_p^\top \b \xi_i + d_p - \b a_p^\top \b x}{\| \b b_p \|_*} \right) \geq s-r_i, &&\hspace{-50pt} \forall i\in [N], p\in [P].  \label{eq:LP1_last}
    \end{align}
\end{subequations}
\cite{chen2023approximations} showed that LA~\eqref{eq:LP_constr_1} represents a family of the ``best'' inner-approximation in the $(\b x, s, \b r)$-feasible region of~\eqref{eq:wdrjcc_abstract} and LA~\eqref{eq:LP_constr_1} becomes exact by optimizing $\b \kappa \coloneqq [\kappa_i]_{i\in[N]}$. 
Furthermore, as demonstrated by~\cite{chen2023approximations}, LA~\eqref{eq:LP_constr_1} with $\b \kappa = \b 1$ is equivalent to W-CVaR for a specific setting of \(\b w\), and becomes exact under certain conditions. These conditions are also extended to our proposed approximation as will be illustrated in Corollaries~\ref{cor:wcvar} and~\ref{cor:exact_conditions}. The $\b x$-feasible region of the LA approximation can be expressed as:
$
    \mathcal{X}_\text{LA}(\b \kappa) = 
    \left\{
          \b x \in \mathcal{X} \mid \exists s, \b r :~\eqref{eq:LP1_s_r}\mbox{--}\eqref{eq:LP1_last}
    \right\}.
$

\section{Strengthened and Faster Linear Approximation}\label{sec:main_theo_res}


\begin{figure}
    \centering
    \begin{tikzpicture}[ node distance=1.2cm and 1.5cm, font=\scriptsize,
                          every node/.style={draw, text width=3cm, align=center, rounded corners},
                          arrow/.style={-{Stealth}}
                        ]

        \node (rhs) {RHS-WDRJCC with Distance Function~\eqref{eq:wdrjcc_abstract}~\citep{2022_OR_exact_DRO}};
        \node (strengthenedR) [right=of rhs] {Strengthened Reformulation with Distance Function~\eqref{eq:MILP_constr_reduced}};
        \node (milp) [above=of strengthenedR] {MIP Reformulation~\eqref{eq:MILP_constr1}~\citep{2022_OR_exact_DRO}};
        \node (strengthenedM) [right=of milp] {Strengthened MIP Reformulation (ExactS)~\eqref{eq:ExactS}~\citep{exact_milp_strengthened}};
        \node (sfla) [below=of strengthenedM, yshift=-0.15cm] {The Proposed SFLA~\eqref{eq:SFLA}};
        \node (la) [below=of sfla] {LA~\eqref{eq:LP_constr_1}~\citep{chen2023approximations}};
        \node (wcvar) [below=of la] {W-CVaR~\eqref{eq:CVaR_RHS}~\citep{mohajerin2018data,chen2023approximations}};

        \draw[arrow] (rhs) to node[midway, above, draw=none] {{\Large $=$}} (strengthenedR);
        \path[arrow] (rhs) to node[midway, below, draw=none] {Prop.~\ref{proposition}} (strengthenedR);
        \draw[arrow, bend left=23] (rhs) to node[midway, above, draw=none] {{\Large $=$}} (milp);
        \draw[arrow, bend right=11] (rhs) to node[midway, below, draw=none] {\Large$\supseteq$} (la);
        \draw[arrow, bend right=22] (rhs) to node[midway, below, draw=none] {\Large$\supseteq$} (wcvar);
        \draw[arrow] (strengthenedR) to node[midway, right, draw=none] { } (sfla);
        \draw[arrow] (milp) to node[midway, above, draw=none] {\Large$=$} (strengthenedM);

        \node at (10.1,-2.6) [draw = none] {\rotatebox{90}{\Large$=$}};
        \node at (10.1,-0.9) [draw = none] {\rotatebox{90}{\Large$=$}};
        \node at (10.1,0.95) [draw = none] {\rotatebox{90}{\Large$=$}};
        \path[arrow] (wcvar) to node[right, draw=none] {$\kappa = \b 1$ and ${\b w}^*$ defined in Cor.~\ref{cor:wcvar}} (la);
        \path[arrow] (sfla) to node[right, draw=none] {$\b \kappa = \b 1$ Cor.~\ref{cor:eq1}} (la);
        \path[arrow] (sfla) to node[right, draw=none] {Conditions in Cor.~\ref{cor:obj_dependent_exact} and~\ref{cor:exact_conditions}} (strengthenedM);

        \node at (8.9,-2.6) [draw = none] {\rotatebox{90}{\Large$\subseteq$}};
        \node at (8.9,-0.9) [draw = none] {\rotatebox{90}{\Large$\subseteq$}};
        \node at (8.9,0.95) [draw = none] {\rotatebox{90}{\Large$\subseteq$}};
        \path[arrow] (wcvar) to node[left, draw=none] {$\kappa = \b 1$ \\and ${\b w} = \frac{1}{P}$ \\Cor.~\ref{cor:wcvar_la}} (la);
        \path[arrow] (sfla) to node[left, draw=none] {Thm.~\ref{theo:propLP_sup_LPori_subMILP}} (la);
        \path[arrow] (sfla) to node[left, draw=none] {Thm.~\ref{theo:propLP_sub_MILP}} (strengthenedM);

    \end{tikzpicture}
    \caption{The derivation flow and the formulation comparisons of the $\b x$-feasible region defined by different reformulation of RHS-WDRJCC.}
\label{fig:flowchart}
\end{figure}

Our proposed approximation starts from the strengthening process of the exact MIP reformulation. As discussed, the RHS-WDRJCC reformulation~\eqref{eq:wdrjcc_abstract} can be exactly reformulated as a set of MIP constraints in~\eqref{eq:MILP_constr1}. \cite{exact_milp_strengthened} showed that the \textcolor{black}{exact MIP reformulation}~\eqref{eq:MILP_constr1} can be further strengthened to another exact MIP reformulation (denoted by ExactS), by replacing $NP-\lfloor \epsilon N \rfloor P$ constraints with $P$ valid inequalities. For convenience, the formulation of ExactS is reproduced as Eq.~\eqref{eq:ExactS} in Appendix~\ref{appen:exactS}. The $\b x$-feasible region of ExactS remains equivalent to the exact MIP set~\eqref{eq:MILP_constr1}, but the number of constraints and the complexity of the whole $(\b x, s, \b r)$-feasible region are reduced. In this section, we will show that this strengthening process can be migrated to the RHS-WDRJCC~\eqref{eq:wdrjcc_abstract} expressed by the distance function (Proposition~\ref{proposition}). \textcolor{black}{Then, combining with the derivation process for LA, we arrive at} our proposed approximation, termed as \textit{Strengthened and Faster Linear Approximation} (SFLA). \textcolor{black}{By combining ExactS and LA}, rather than compromising the $\b x$-feasible region as illustrated in Case I of Figure~\ref{fig:valid_ineq}, we prove that the proposed SFLA is always the Case II: the proposed SFLA will not introduce additional conservativeness and can even be less conservative than LA proposed by \cite{chen2023approximations} (Theorem \ref{theo:propLP_sup_LPori_subMILP}). 
\textcolor{black}{The proof of all the theoretical results in this section is provided in Appendix~\ref{appen:proof}.}

{\color{black}
\subsection{Derivation and Formulation}

Let $k \coloneqq \lfloor \epsilon N \rfloor$ and assume an ordering $\b b_p^\top \b \xi_1 \leq \cdots \leq \b b_p^\top \b \xi_N$ without loss of generality. We further denote by $q_p$ the $(k+1)$-th smallest value as $q_p \coloneqq \b b_p^\top \b \xi_{k+1}$.
Then we introduce the index set $[N]_p$ with only $k = \lfloor \epsilon N \rfloor$ elements as $[N]_p \coloneqq \left\{ i\in [N] \mid  \b b_p^\top \b \xi_i < q_p \right\}$ for $p\in [P]$. We propose the following proposition which demonstrates that the strengthening process for the exact MIP reformulation ExactS~\eqref{eq:MILP_constr1} in~\cite{exact_milp_strengthened} can be migrated to the exact RHS-WDRJCC reformulation~\eqref{eq:wdrjcc_abstract} expressed by the distance function $\mathrm{dist}(\b \xi_i, \mathcal{S}(\b x))$.
\begin{proposition}\label{proposition}
The feasible set $\mathcal{X}_\text{Exact}$ defined through constraints \eqref{eq:wdrjcc_abstract} can be equivalently defined by keeping only $\lfloor \epsilon N \rfloor P$ constraints in~\eqref{eq:dist_constr} as~\eqref{eq:MILP_reduced_part} and replacing the remaining $NP-\lfloor \epsilon N \rfloor P$ constraints in~\eqref{eq:dist_constr} involving the distance function with $P$ valid inequalities~\eqref{eq:MILP_reduced_q}:
\begin{subequations}\label{eq:MILP_constr_reduced}
    \begin{align}
        & \left(\dfrac{\b b_p^\top \b \xi_i + d_p - \b a_p^\top \b x}{\| \b b_p \|_*} \right)^+  \geq s - r_i, &&\hspace{-50pt} \forall i \in [N]_p, p\in [P], \label{eq:MILP_reduced_part}\\
        & \dfrac{q_p + d_p - \b a_p^\top \b x}{\| \b b_p \|_*} \geq s, &&\hspace{-50pt} \forall p\in [P].\label{eq:MILP_reduced_q}
    \end{align}
\end{subequations}
In other words,
$
    \mathcal{X}_{\text{Exact}} =
        \{
        \b x \in \mathcal{X} \mid \exists s, \b r :~\eqref{eq:s_r}, \eqref{eq:epsilonNs}, \eqref{eq:MILP_reduced_part}, \eqref{eq:MILP_reduced_q}
        \} = 
    \{
        \b x \in \mathcal{X} \mid 
        \exists s, \b r :~\eqref{eq:s_r}\mbox{--}\eqref{eq:dist_constr}
    \}
$.
\end{proposition}

\noindent The structure of the constraints~\eqref{eq:MILP_constr_reduced} motivates SFLA, which linearizes the non-convex constraints in~\eqref{eq:MILP_reduced_part} using the approximate distance function $\widehat{\mathrm{dist}}(\b \xi_i, \mathcal{S}(\b x))$~\eqref{eq:dist_approx} as was done for LA, while preserving the strengthened formulation in~\eqref{eq:MILP_constr_reduced}. Therefore, SFLA can be expressed as follows: 
\begin{subequations}\label{eq:SFLA} 
    \begin{align}
        & s \geq 0, \b r \geq \b 0, \label{eq:prop_LP_s_r} \\
        & \epsilon N s - \sum\limits_{i\in [N]}r_i \geq \theta N, \label{eq:prop_LP_epsilonNs} \\
        & \kappa_i \left(\dfrac{\b b_p^\top \b \xi_i + d_p - \b a_p^\top \b x}{\| \b b_p \|_*} \right) \geq s - r_i, &&\hspace{-50pt} \forall i \in [N]_p,\ p\in [P], \label{eq:prop_LP_part}\\
        & \dfrac{q_p + d_p - \b a_p^\top \b x}{\| \b b_p \|_*} \geq s, &&\hspace{-50pt} \forall p\in [P], \label{eq:prop_LP_q}
    \end{align}
\end{subequations}
where we still have the preset hyperparameters $\kappa_i \in [0, 1]$ for all $i \in [N]$. The $\b x$-feasible region of the proposed SFLA can be expressed as $\mathcal{X}_\text{SFLA}(\b\kappa) \coloneqq \left\{ \b x \in \mathcal{X} \mid  \exists s, \b r :~\eqref{eq:prop_LP_s_r}\mbox{--}\eqref{eq:prop_LP_q} \right\}$.

\subsection{Key Properties and Practical Relevance}
\label{subsec:key_properties}

By integrating the strengthening process used in ExactS~\citep{exact_milp_strengthened} with the approximate distance function in LA~\citep{chen2023approximations}, SFLA is not only more computationally efficient but also less (or equally) conservative than existing convex approximation approaches, making it well-suited for practical large-scale settings.

\subsubsection{Computational Tractability}
\label{subsubsec:tractability}

Many real-world decision-making problems prioritize computational tractability over strict optimality due to limited time budgets (while still requiring good solution quality). As discussed in Section~\ref{sec:intro}, this is the case for applications such as unit commitment and multilevel optimization with chance-constrained lower levels, where even the deterministic models are already computationally challenging, and convexity is important for tractable reformulation. In these settings, the value of approximation of chance constraints hinges not only on its probabilistic guarantees but also on whether it can be embedded into a large optimization pipeline and solved reliably at scale.

The proposed SFLA is therefore of high practical importance for two reasons. First, compared with the exact MIP reformulation ExactS, SFLA \emph{convexifies} and \emph{linearizes} the WDRJCC~\eqref{constr:wdrjcc} through an approximate distance function. This convexity is crucial for maintaining reformulation tractability, scalability, and numerical robustness, especially when SFLA is applied to the lower-level JCC of multi-level problems or repeatedly solved within a decomposition and rolling-horizon procedure.

Second, compared with common convex approximation benchmarks such as LA~\eqref{eq:LP_constr_1} and W-CVaR~\eqref{eq:CVaR_RHS} for WDRJCC, SFLA admits a \emph{more compact} representation. Due to the migration of the strengthening process, the valid inequality \eqref{eq:prop_LP_q} cuts off redundant space of ancillary variables and thus brings computational speedup. In addition, this valid inequality enables replacing the original $[N]$ index set with $[N]_p$ in constraint~\eqref{eq:prop_LP_part}. Due to this replacement, the proposed SFLA only has $\lfloor \epsilon N \rfloor P + P +1$ constraints (except for basic domain constraints such as non-negativity), while WCVaR has $NP+P+1$ constraints and LA has $NP+1$ constraints. For decision-making problems that require high reliability, $\epsilon$ will be set small, which reduces the number of constraints significantly. This will, in turn, reduce the memory footprint and is important for large-scale problems.

\subsubsection{Safety and No Extra Conservativeness}
\label{subsubsec:safe_not_conservative}

Although computational tractability is prioritized, a good (although not necessarily optimal) solution quality is still important. Also, as an approximation, real-world applications that demand high reliability often prefer inner-approximation over outer-approximation due to safety. This safety property is preserved for the proposed SFLA:
\begin{theorem}
\label{theo:propLP_sub_MILP}
    The set $\mathcal{X}_\text{SFLA}(\b \kappa)$ is a subset of $\mathcal{X}_\text{Exact}$, i.e., $\mathcal{X}_\text{SFLA}(\b \kappa) \subseteq \mathcal{X}_\text{Exact}$.
\end{theorem}

Regarding solution quality, we show that although being more computationally tractable than the exact MIP reformulation and the two common convex approximations, including LA and W-CVaR, the proposed SFLA will not bring extra conservativeness than LA and W-CVaR, and can even be exact for certain conditions. 

\noindent First, we show that SFLA does not introduce additional conservativeness compared to LA:
\begin{theorem}
\label{theo:propLP_sup_LPori_subMILP}
    The set $\mathcal{X}_\text{SFLA}(\b\kappa)$ is a superset of $\mathcal{X}_\text{LA}(\b\kappa)$, i.e., $\mathcal{X}_\text{LA}(\b\kappa) \subseteq \mathcal{X}_\text{SFLA}(\b\kappa)$.
\end{theorem}

Moreover, under a natural choice of the parameter, the two formulations coincide:
\begin{corollary}
\label{cor:eq1}
    We have $\mathcal{X}_\text{SFLA}(\b 1) = \mathcal{X}_\text{LA}(\b 1)$.
\end{corollary}
Second, SFLA with $\b\kappa=\b 1$ dominates W-CVaR with $\b w=\b 1/P$, which is the standard W-CVaR in~\citet{chen2023approximations} and~\citet{jiang2024also}:
\begin{corollary}\label{cor:wcvar_la}
We have $\mathcal{X}_\text{WCVaR}(\b 1 / P) \subseteq \mathcal{X}_\text{SFLA}(\b 1)$.
\end{corollary}

In addition, they coincide when W-CVaR has tuned hyperparameters:
\begin{corollary}\label{cor:wcvar}
    We have $\mathcal{X}_\text{SFLA}(\b 1) = \mathcal{X}_\text{WCVaR}(\b w^*)$, where $\b w^*=[w^*_p]_{p\in [P]}$ and $w^*_p = \frac{\| \b b_p \|_*^{-1}}{\sum_{l \in [P]} \| \b b_l \|_*^{-1}}$.
\end{corollary}
Finally, SFLA can recover the exact feasible set by allowing $\b\kappa$ to be decision-dependent:
\begin{corollary}\label{cor:obj_dependent_exact}
    Let $\b \kappa \in [0,1]^N$ be extra decision variables and define the updated SFLA set as $\mathcal{X}_\text{SFLA}^* = \left\{ \b x \in \mathcal{X} \mid  \exists \b \kappa \in [0,1]^N, s, \b r:~\eqref{eq:prop_LP_s_r}\mbox{--}\eqref{eq:prop_LP_q} \right\}$.
    Then we have $\mathcal{X}_\text{SFLA}^* = \mathcal{X}_\text{Exact}$.
\end{corollary}
This exactness primarily serves as a theoretical benchmark for understanding the tightness of SFLA since treating $\b\kappa$ as extra decision variables renders the problem nonconvex and bilinear. Nevertheless, even with fixed $\b\kappa=\b 1$, exactness can still arise:
\begin{corollary}\label{cor:exact_conditions}
    The equality $\mathcal{X}_\text{SFLA}(\b 1) = \mathcal{X}_\text{Exact}$ holds under either of the following conditions: (i) We have $\b\xi_i \in \mathcal{S}(\b x)$ for all $i \in [N]$ and $\b x \in \mathcal{X}_\text{Exact}$; (ii) We have $\epsilon \leq 1/N$.
\end{corollary}
The second condition in Corollary~\ref{cor:exact_conditions} shows that the proposed SFLA becomes exact when the number of data $N$ is small or the risk level $\epsilon$ is low. This also implies that a smaller $N$ or a smaller $\epsilon$ can lead to better approximation quality for the proposed SFLA. This implication is particularly important, as DRCCP is often used in cases with limited historical data due to its strong out-of-sample performance. Additionally, safety-critical applications, such as in power systems, generally demand a low risk level. Thus, the proposed SFLA can achieve high approximation quality across a range of practical applications.

Figure~\ref{fig:flowchart} summaries the set relationships between the proposed SFLA and other reformulations. These many theoretical properties regarding conservativeness can hold for the proposed SFLA because SFLA combines the derivation process of the both exact reformulation~\citep{exact_milp_strengthened} and LA~\citep{chen2023approximations}. A detailed proof can be found in Appendix \ref{appen:proof}.

Finally, we provide a simple prior guideline for setting the hyperparameters. The following monotonicity property implies that $\b\kappa=\b 1$ yields the largest (least conservative) inner approximation within this family, and hence serves as a natural default choice when no additional tuning information is available:
\begin{corollary}\label{cor:hyper}
    We have $\mathcal{X}_{SFLA}(\b 1) \supseteq \mathcal{X}_{SFLA}(\kappa \b 1)$ for all $\kappa \in [0, 1]$.
\end{corollary}

}

{\color{black}
\section{Robustness Maximization}\label{sec:robustsmax}

In WDRJCC, the risk level $\epsilon$ and the Wasserstein radius $\theta$ do admit clear interpretations in terms of risk and ambiguity. However, for decision-makers that demand high reliability but have other backgrounds (e.g., power system operators), it can be easier and more interpretable to first specify a bearable cost as a constraint that bounds the original objective, and then minimize the risk level $\epsilon$ or maximize the Wasserstein radius $\theta$ to align with their high reliability principle. This paper refers to this alternative decision-making paradigm as \emph{robustness maximization}. In this section, we illustrate how the proposed SFLA can be applied within this paradigm and clarify the link between minimizing the risk level $\epsilon$ and maximizing the Wasserstein radius $\theta$.

Note that our robustness maximization framework can be related to the concept of \emph{robust satisficing} \citep{schwartz2011makes}, in which a decision-maker specifies an acceptable utility target and then seeks to maximize the robustness of achieving it. A direct formulation of this concept treats the risk level $\epsilon$ as an additional decision variable and minimizes $\epsilon$ in the objective \citep[Online Supplement~D]{zhao2023distributionally}. \citet{long2023robust} further formalize this concept by proposing a framework in which the \emph{adjusted} bearable cost is required to hold for \emph{all} distributions, with the adjustment determined by their distance from a reference distribution.

}

\subsection{Risk Minimization}

In this formulation, the objective is to minimize the risk level $\epsilon$:
\begin{subequations}
\label{eq:eps_robust_satisfi}
\begin{align}
    \min_{\b x \in \mathcal{X}, 0\leq \epsilon \leq \epsilon_{\max}} \quad & \epsilon \label{obj:eps_robust_satisfi}  \\
    \mbox{s.t.} \quad & \sup_{\mathbb {P} \in \mathcal{F}_N (\theta)} \mathbb{P}[\b \xi \notin \mathcal{S}(\b x)] \leq \epsilon , \label{constr:wdrjcc_eps_robust_satisfi} \\
    & c(\b x) \leq \tau \label{constr:eps_robust_satisfi},
\end{align}
\end{subequations}
where $\epsilon_{\max} \in [0, 1]$ is an upper bound on the risk level $\epsilon$, determined by decision maker's risk tolerance. Then constraint~\eqref{constr:eps_robust_satisfi} ensures that the cost $c(\b x)$ does not exceed the cost level $\tau$. Moreover, choosing a smaller $\epsilon_{\max}$ can also tighten the problem to speedup the computation. \textcolor{black}{Note that, because our objective is to maximize robustness (i.e., to minimize $\epsilon$), the bearable cost constraint~\eqref{constr:eps_robust_satisfi} already regulates the optimization. Therefore, there is no need to explicitly penalize $\epsilon$ as the case for \citet{zhang2025integer}---where the risk level is intended to increase.}

The Exact reformulation \eqref{eq:wdrjcc_abstract}, WCVaR \eqref{eq:CVaR_RHS}, and LA \eqref{eq:LP_constr_1} can be easily extended to this setting, by taking $\epsilon \in [0, \epsilon_{\max}]$ as an extra decision variable and making corresponding constraints \eqref{eq:epsilonNs}, \eqref{eq:LP1_epsilonNs} and \eqref{eq:wcvar-b} bilinear and nonconvex. Take LA as an example, the $\b x$-feasible region of the extended LA is denoted as $\mathcal{X}_\text{LA}'(\b\kappa) \coloneqq \{ \b x \in \mathcal{X} \mid \exists s, \b r, \epsilon:~\eqref{eq:LP1_s_r}\mbox{--}\eqref{eq:LP1_last}, \epsilon \in [0, \epsilon_{\max}]\}$.

The risk minimization formulation is not directly applicable to the proposed SFLA in~\eqref{eq:SFLA} and ExactS in~\eqref{eq:ExactS}, because the valid inequalities and the index set $[N]_p$ for these two sets depend on $\epsilon$ (recall the ranking-related definition where $k \coloneqq \lfloor \epsilon N \rfloor$), which, if formulated, leads to complicated non-convex combinatorial formulations. However, given $\epsilon_{\max}$, we propose to extend SFLA that works for the risk minimization problem, by replacing constraints \eqref{eq:prop_LP_part} and \eqref{eq:prop_LP_q} with the following:
\vspace{-7mm}
\begin{subequations}\label{eq:SFLA_eps}
    \begin{align}
        & \kappa_i \left(\dfrac{\b b_p^\top \b \xi_i + d_p - \b a_p^\top \b x}{\| \b b_p \|_*} \right) \geq s - r_i, &&\hspace{-50pt} \forall i \in [N]_p',\ p\in [P], \label{eq:prop_LP_part_eps}\\
        & \dfrac{q_p' + d_p - \b a_p^\top \b x}{\| \b b_p \|_*} \geq s, &&\hspace{-50pt} \forall p\in [P]. \label{eq:prop_LP_q_eps}
    \end{align}
\end{subequations}
where we define $q_p' \coloneqq \b b_p^\top \b \xi_{k'+1}$ and $[N]_p' \coloneqq \left\{ i\in [N] \mid  \b b_p^\top \b \xi_i < q_p' \right\}$, with $k' \coloneqq \lfloor \epsilon_{\max} N \rfloor$ and the same ordering assumption $\b b_p^\top \b \xi_1 \leq \cdots \leq \b b_p^\top \b \xi_N$ without loss of generality. In this way, we eliminate the dependency on $\epsilon$. We further denotes its $\b x$-feasible region as $\mathcal{X}_\text{SFLA}'(\b\kappa) \coloneqq \{ \b x \in \mathcal{X} \mid \exists s, \b r, \epsilon:~\eqref{eq:prop_LP_s_r},\eqref{eq:prop_LP_epsilonNs},\eqref{eq:prop_LP_part_eps},\eqref{eq:prop_LP_q_eps}, \epsilon \in [0, \epsilon_{\max}] \}$.
We have the following result:
\begin{corollary}\label{theo:min_eps}
    $\mathcal{X}_\text{LA}'(\b\kappa) \subseteq \mathcal{X}_\text{SFLA}'(\b\kappa)$.
\end{corollary}
The proof can follow the same process as the proof for Theorem \ref{theo:propLP_sup_LPori_subMILP}, noticing the fact that $q'_p \geq q_p$ and $[N]_p \subseteq [N]_p'$, such that the derivation process that holds for constraints \eqref{eq:prop_LP_part} and \eqref{eq:prop_LP_q} in $\mathcal{X}_\text{SFLA}(\b\kappa)$ will also hold for \eqref{eq:prop_LP_part_eps} and \eqref{eq:prop_LP_q_eps} in $\mathcal{X}_\text{SFLA}'(\b\kappa)$.

\subsection{Radius Maximization}\label{sec:radius_maximization}

To maximise robustness, one may alternatively maximize the radius of the ambiguity set, subject to the same tolerable cost constraint~\eqref{constr:theta_robust_satisfi}:
\begin{subequations}
\label{eq:theta_robust_satisfi}
\begin{align}
    \max_{\b x \in \mathcal{X}, \, \theta_{\min} \leq \theta} \quad & \theta \label{obj:theta_robust_satisfi}  \\
    \text{s.t.} \quad & \sup_{\mathbb{P} \in \mathcal{F}_N(\theta)} \mathbb{P}[\b \xi \notin \mathcal{S}(\b x)] \leq \epsilon , \label{constr:wdrjcc_theta_robust_satisfi} \\
    & c(\b x) \leq \tau. \label{constr:theta_robust_satisfi}
\end{align}
\end{subequations}
This formulation is closely related to the risk minimization formulation \eqref{eq:eps_robust_satisfi}. Specifically, let $\mathcal{X}_\text{risk}(\theta)$ denote the $\b x$-feasible region of \eqref{eq:eps_robust_satisfi} for a fixed $\theta$, and let $\mathcal{X}_\text{rad}(\epsilon)$ denote the $\b x$-feasible region of \eqref{eq:theta_robust_satisfi} for a fixed $\epsilon$. Then we have:
\begin{theorem}\label{theo:equivalence_risk_radius}
    Given $\theta_{\min} \geq 0$ and $\epsilon_{\max} \in [0, 1]$, we have
    $\mathcal{X}_\text{risk}(\theta_{\min}) = \mathcal{X}_\text{rad}(\epsilon_{\max})$.
\end{theorem}

The proof only focuses on the constraint $\epsilon N s - \sum\limits_{i \in [N]} r_i \geq \theta N$, which appears in most WDRJCC reformulations such as the proposed SFLA, the exact reformulation, ExactS, LA, and WCVaR. We can show that any feasible solution $(\b x, s, \b r)$ to one formulation is also feasible to the other under appropriately chosen values of $\theta$ and $\epsilon$. Therefore, Theorem \ref{theo:equivalence_risk_radius} is not restricted to the exact reformulation but extends naturally to SFLA, LA, WCVaR, and ExactS.

Theorem \ref{theo:equivalence_risk_radius} has several implications. On the one hand, although treating $0\leq \epsilon \leq \epsilon_{\max}$ as an extra decision variable will make the problem bilinear and non-convex, the $\b x$-feasible space (or the $(\b x, s, \b r)$-feasible space for reformulations such as SFLA) is still convex and linear (because treating $\theta$ as an extra decision variable is still linear and convex). This may promote a broader range of optimization techniques for solving problems with $\epsilon$ being decision variables, for example, still minimizing the original cost function but plus a penalty cost for the increased risk \citep{zhang2025integer}, or possibly enlightening the methods for tuning the individual risk levels $\epsilon_p$ in the Bonferroni approximation (this is more complicated due to the coupling $\sum_{p\in[P]}\epsilon_p\leq \epsilon$) \citep{ding2022distributionally}.

On the other hand, Theorem~\ref{theo:equivalence_risk_radius} implies that the risk minimization model~\eqref{eq:eps_robust_satisfi} and the radius maximization model~\eqref{eq:theta_robust_satisfi} share the same $x$-feasible region. Therefore, comparing the two formulations reduces to evaluating which objective---minimizing $\epsilon$ or maximizing $\theta$---is a more \emph{effective} criterion for selecting a \emph{robust} solution.  First, radius maximization is generally more computationally tractable, as it preserves the linearity and convexity of the constraint $\epsilon N s - \sum_{i \in [N]} r_i \geq \theta N$. This allows reformulation techniques such as ExactS, the proposed SFLA, LA, and WCVaR to be directly applied by treating $\theta$ as a one-dimensional decision variable. Second, our numerical results (see Appendix~\ref{res:risk_vs_radius}) show that radius maximization tends to achieve slightly higher out-of-sample robustness than risk minimization. Therefore, for robustness maximization, we advocate using radius maximization over risk minimization.

\section{Numerical Experiments}
\label{sec:casestudies}

This section presents numerical case studies to demonstrate the reduction of computational complexity and also the approximation quality of the proposed SFLA. For all norm calculations in the RHS-WDRJCC formulations, we adopt the L2 norm which has the property of being self-dual. We set $\b \kappa=\b 1$ for the proposed SFLA and LA, $\b w_p=1/P$ for W-CVaR following \cite{ordoudis2021energy}, and $\epsilon_p = \epsilon / P$ for Bonferroni approximation. To reach robust conclusions, we implement multiple random runs for each parameter setting as detailed in Appendices~\ref{appen:setting_UC} and \ref{appen:setting_bilevel}. The code, the mathematical formulations, and detailed case study settings for our UC and bilevel problems are available at our GitHub repository \citep{Zhou2024_code}. 

\subsection{Unit Commitment Problem}
\label{sec:UC}



The first case study examines a chance-constrained UC problem. UC represents a typical real-world large-scale problem, and its primary objective is to determine the generator on/off statuses for the next scheduling period with a minimal operation cost~\citep{MISO_UC_difficulty}. UC is typically formulated as an MIP problem, where each generator requires at least $T$ binary variables to represent on/off statuses over the next $T$ time steps. With thousands of generators (e.g., $1{,}400$ in MISO), non-linear generator cost functions~\citep{xu2017application}, and complex network constraints, UC is challenging to solve even in a deterministic case. However, the integration of renewables introduces additional uncertainty into UC, making chance constraints essential to ensure reliability~\citep{wang2016risk, wu2016solution, yang2019analytical}. Given UC's inherent complexity, minimizing the added computational burden of JCCs is critical.

We compare the proposed SFLA with benchmarks including the exact and strengthened MIP formulation (ExactS)~\citep{exact_milp_strengthened} in Eq.~\eqref{eq:ExactS}, LA~\citep{chen2023approximations} in Eq.~\eqref{eq:LP_constr_1}, W-CVaR~\citep{mohajerin2018data, chen2023approximations} in Eq.~\eqref{eq:CVaR_RHS}, Bonferroni approximation~\citep{chen2023approximations} in Eq.~\eqref{eq:Bonfe_RHS}, and ALSO-X\# \citep{jiang2024also}. 

\subsubsection{Evaluation Metrics}
\label{sec:uc_metrics}

The UC problem is primarily evaluated using the following metrics:

\noindent 1) \emph{TimeF} (s): This metric records the time (in seconds) to obtain the first ``comparable high-quality'' solution, defined as the solution with an objective within the \texttt{MIPGap} of the final objective achieved by the proposed SFLA (whether optimal or the best found within the \texttt{TimeLimit}). If SFLA does not yield a feasible solution within the \texttt{TimeLimit}, the \emph{TimeF} of the proposed SFLA is set to the \texttt{TimeLimit}, and the \emph{TimeF} for other benchmarks is set to the time when the first feasible solution is found (if any, otherwise set to the \texttt{TimeLimit}). Because MIP admits early stop, this metric offers a practical perspective in cases where the solver identifies a solution close to optimal early but requires extra time for verifying optimality.

\noindent 2) \emph{Time} (s): This measures the total time (in seconds) spent to solve the problem to optimality (within \texttt{MIPGap}), and is set to the \texttt{TimeLimit} if no optimal solution is found within the \texttt{TimeLimit}. This metric compliments the previous \emph{TimeF} by providing verified optimality.

\noindent 3) \emph{Obj. Diff.} (\%): This is the percentage difference between the optimal value (minimum UC cost) of a benchmark method and that of the SFLA, relative to the optimal value of the SFLA. A positive value indicates that the benchmark achieves a lower cost (better optimality) than the SFLA. We only keep certain random runs where the benchmark method and the proposed SFLA are both solved to optimality. 

\noindent 4) \emph{Reli.} (\%): This is the out-of-sample joint satisfaction rate of the RHS-WDRJCC. This metric is used to justify the reasonability of the selection of Wasserstein radius $\theta$, and is calculated on $5000$ held-out samples.

This section prioritizes \emph{TimeF} over the total computation time \emph{Time} due to comparison with ExactS. MIP can always employ early stopping, and ExactS can achieve better optimality compared to inner approximations (SFLA, LA, W-CVaR, Bonferroni, and ALSO-X\#). Therefore, a higher \emph{Time} for ExactS does not necessarily suggest its disadvantage.

\subsubsection{Results for Unit Commitment}

\begin{figure}
\vspace{-4mm}
     \FIGURE
    {\includegraphics[width=1\linewidth]{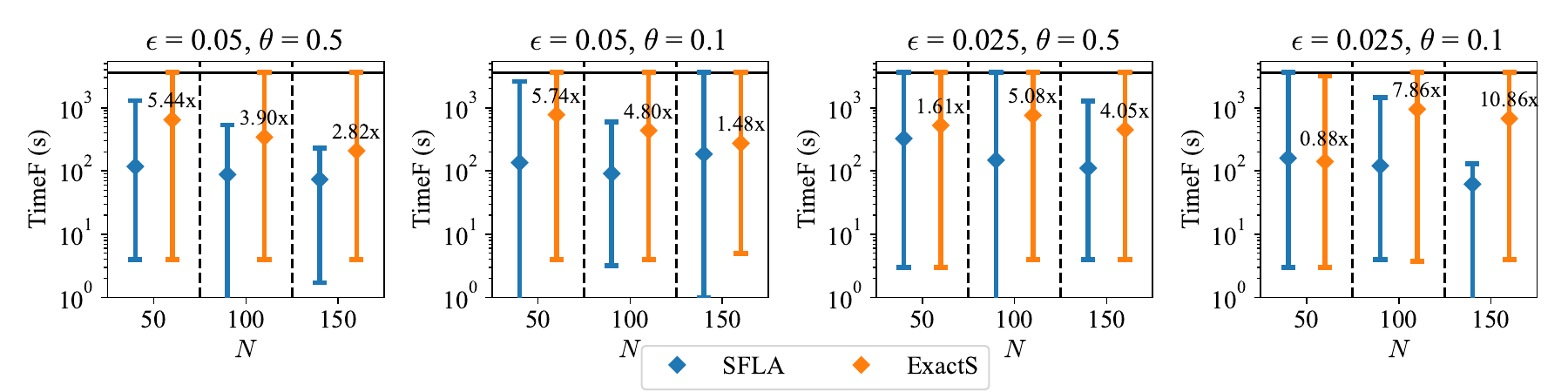}\vspace{-2mm}} 
    {Computation time to obtain the first comparable high-quality solution for the proposed SFLA and the benchmark ExactS for the UC problem. \label{fig:suc_timeF}}
    {Dots represent the mean value of the 150 random runs, with error bars indicating the 95\% percentile interval (from 2.5th to 97.5th). The black horizontal line represents the 3600s \texttt{TimeLimit}. The optimization horizon is set to $T = 24$. Numbers ending with ``x'' represent the speedup in the average computing time of the proposed SFLA compared to ExactS.}
\end{figure}
\begin{figure}
     \FIGURE
    {\includegraphics[width=1\linewidth]{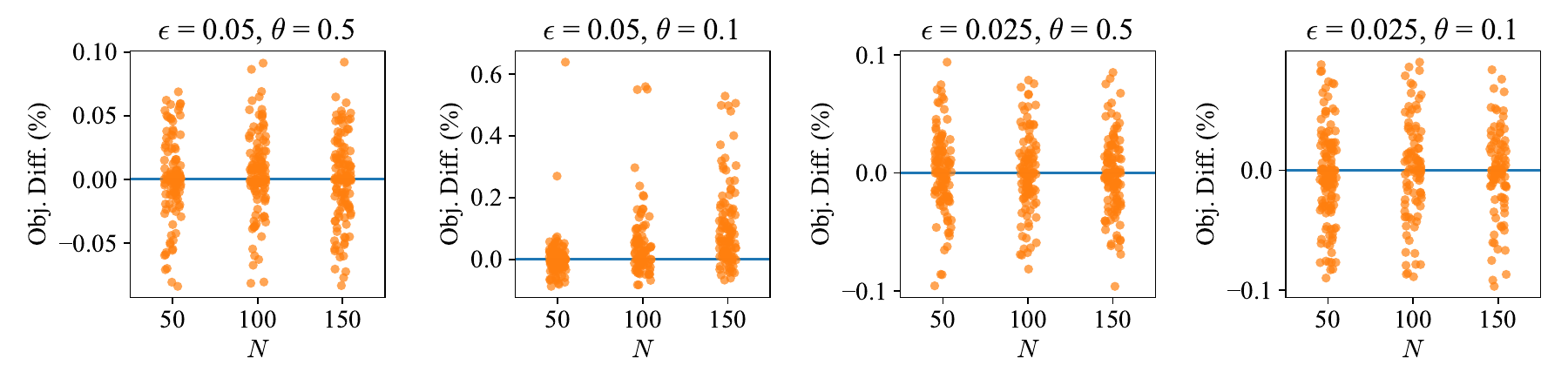} \vspace{-5mm}} 
    {Performance comparison of the optimality for the proposed SFLA and the benchmark ExactS for the UC problem. \label{fig:suc_gap}}
    {Each round dot represents a result of one of the 150 random runs. The blue horizontal line indicates a zero difference of the optimal value of a benchmark method compared to the proposed SFLA. Values higher than this horizontal line indicate that the benchmark achieves \textcolor{black}{a} lower cost (better optimality) than the proposed SFLA. The optimization horizon is set to $T = 24$.}
\end{figure}

Our main content specifically focuses on the comparison between the proposed SFLA and ExactS. Comparisons with LA, W-CVaR, Bonferroni, and ALSO-X\# are excluded here due to their significant disadvantages relative to the proposed SFLA as demonstrated in Appendix~\ref{appen:sup_res_UC}. Also, our numerical results suggest that both the proposed SFLA and ExactS exhibit variability. To reach a more robust conclusion, we test a range of parameter settings and implement $150$ random runs for each setting. Figure~\ref{fig:suc_timeF} visualizes the distribution of the \emph{TimeF} metric. As can be seen, except for a single parameter setting $(\epsilon=0.025, \theta=0.1, N=50)$ where the proposed SFLA shows a minor underperformance ($12\%$ slower), the SFLA demonstrates significant computational speedup in \emph{TimeF} over ExactS: the computational speedup reaches $10\times$ for a parameter setting $(\epsilon=0.025,\ \theta=0.1,\ N=150)$, and the speedup is $3.8\times$ for the average \emph{TimeF} of the 12 parameter settings. 



Figure~\ref{fig:suc_gap} compares the quality of the optimal solution of the proposed SFLA with that of ExactS. In most random runs and most parameter settings, the \emph{Obj. Diff.} is within the preset \texttt{MIPGap}, demonstrating comparable optimal solution quality between the proposed SFLA and ExactS. Slight underperformance (less than 0.65\%) for certain random runs appears only under the combination of a high risk (\(\epsilon=0.05\)), a small radius (\(\theta=0.1\)), and a large number of samples ($N \geq 100$), which is consistent with Corollary \ref{cor:exact_conditions} and the related discussions. However, among these underperformed parameter settings, the average underperformance is still below $0.107\%$. In addition, real-world applications such as those found in power systems have a small tolerance for risk, making the proposed SFLA an effective approximation method.

Our Appendix~\ref{appen:sup_res_UC} further compares the proposed SFLA and ExactS in different problem scales (different $T$). It is observed that the proposed SFLA preserves computational efficiency in more complex cases without sacrificing optimality, making it particularly advantageous for practical high-complexity applications such as UC.

\subsection{Bilevel Storage Strategic Bidding Problem}
\label{sec:bilevel}

This section explores a bilevel optimization problem that allows a price-maker energy storage operator to maximize its profit by participating in a joint day-ahead energy and reserve market. The market participation process can be summarized as follows: Participants submit their bids and offers for each energy and reserve product for specific future delivery times. Bids represent the provision of services to the power system, such as generating power or supplying reserve capacity, while offers involve procuring services, such as purchasing energy for storage charging. In line with most established power markets, we model each bid or offer as a price-quantity pair that is submitted to the market for clearing, which then determines a market-clearing price (marginal price) and cleared quantities bounded by the submitted bid-offer quantities. Upon settlement, market participants receive payment or incur charges at the market-clearing price for their respective cleared quantities, with a binding obligation to fulfill these quantities.

An energy storage operator transitions from a price-taker to a price-maker when its storage power capacity is large enough to influence market outcomes. To maximize profit, the storage operator determines the optimal bidding or offering strategy by solving a bilevel strategic bidding problem, where the lower level models the market-clearing process~\citep{8036231,TOMASSON2020114251,DIMITRIADIS2022123026,7828161}. With increasing uncertainty due to renewable generation, the market clearing process can incorporate chance constraints (RHS-WDRJCC), such as those in the UC problem (see Section~\ref{sec:UC}). Therefore, to accurately model the market-clearing process, the bilevel problem needs to incorporate an RHS-WDRJCC in the lower level. Consequently, convex approximations to RHS-WDRJCC become important for transforming the bilevel problem into a solvable single-level counterpart via KKT conditions or strong duality. Note that chance constraints at lower levels frequently arise in bilevel and multilevel optimization problems, where uncertainty is revealed after lower-level decisions are made~\citep{beck2023survey}. Another example is the European multilevel gas market clearing problem~\citep{heitsch2022convex}.

This section compares the proposed SFLA with typical convex approximation benchmarks, including LA~\citep{chen2023approximations} in Eq.~\eqref{eq:LP_constr_1}, W-CVaR~\citep{mohajerin2018data, chen2023approximations} in Eq.~\eqref{eq:CVaR_RHS}, and Bonferroni approximation~\citep{chen2023approximations} in Eq.~\eqref{eq:Bonfe_RHS}. For all the RHS-WDRJCC approximations, the corresponding bilevel problems are reformulated as single-level mixed-integer nonlinear non-convex programming problems, where the lower-level problem is replaced by the corresponding KKT optimality conditions, which can be solved using the \texttt{Gurobi} solver. ExactS is not included because bilevel problem solving requires convexity at the lower level (ALSO-X\# is also not applicable due to its iterative nature). The bilevel formulation and also the corresponding single-level reformulation using KKT can be found in our online GitHub repository.

\begin{figure}
    \centering    
    \begin{minipage}{\textwidth}
    \vspace{-3.mm}
    \includegraphics[width=\columnwidth]{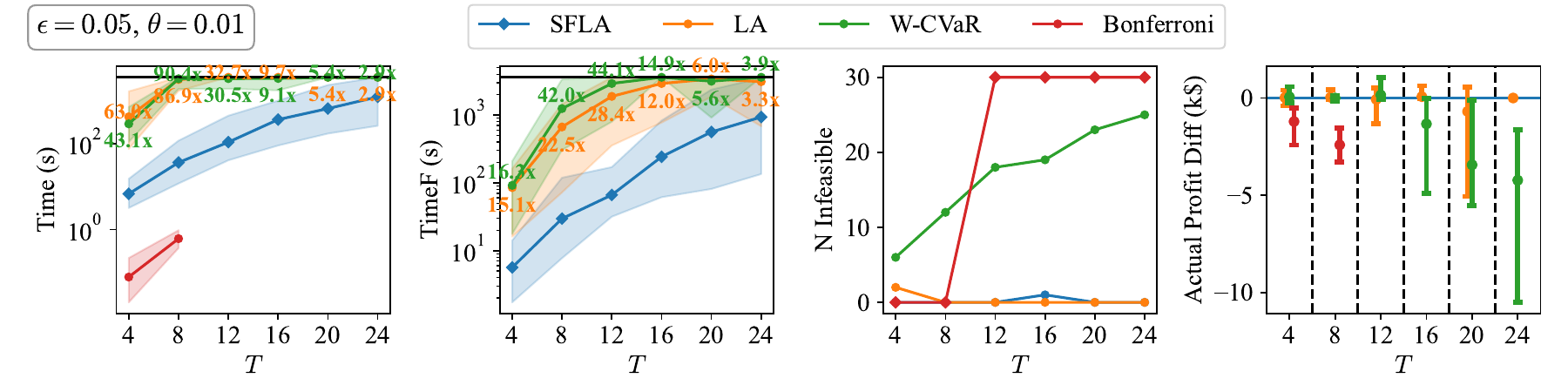}
    \end{minipage}
    \begin{minipage}{\textwidth}
    \hspace{1em}
    \begin{minipage}{\textwidth}
        \centering
        \vspace{-1.6em}  
        \subfloat[\label{fig:time_eps0.05_theta0.01}]{\phantom{\rule{0.247\textwidth}{0pt}}}
        \subfloat[\label{fig:timeF_eps0.05_theta0.01}]{\phantom{\rule{0.247\textwidth}{0pt}}}
        \subfloat[\label{fig:n_infeasible_eps0.05_theta0.01}]{\phantom{\rule{0.247\textwidth}{0pt}}}
        \subfloat[\label{fig:profit_diff_eps0.05_theta0.01}]{\phantom{\rule{0.247\textwidth}{0pt}}}
    \end{minipage}
    \end{minipage}
    \begin{minipage}{\textwidth}
    \vspace{0em}
    \includegraphics[width=\columnwidth]{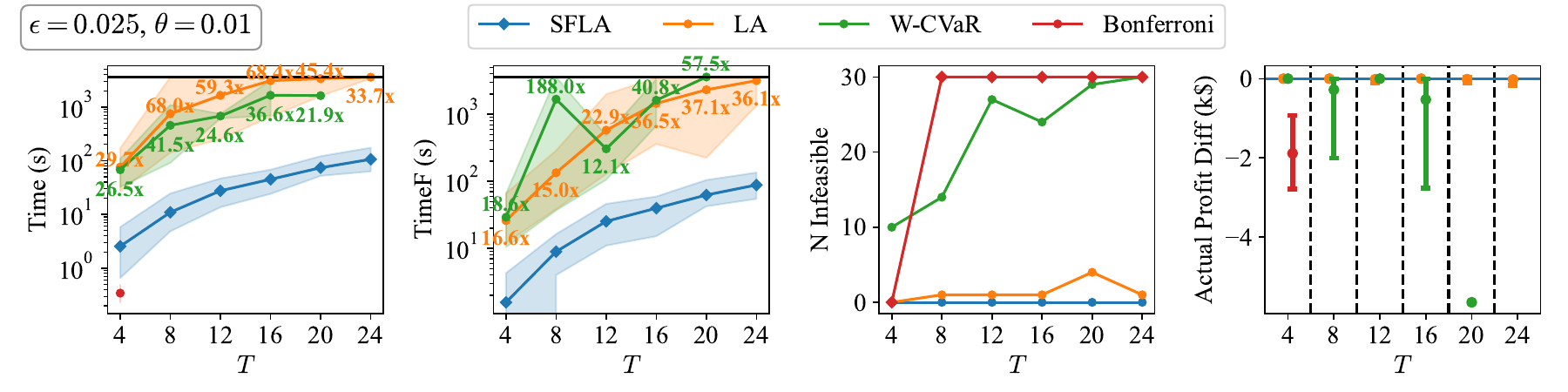}
    \end{minipage}
    \begin{minipage}{\textwidth}
    \hspace{1em}
    \begin{minipage}{\textwidth}
        \centering
        \vspace{-1.5em}  
        \subfloat[]{\phantom{\rule{0.247\textwidth}{0pt}}}
        \subfloat[]{\phantom{\rule{0.247\textwidth}{0pt}}}
        \subfloat[\label{fig:n_infeasible_eps0.025_theta0.01}]{\phantom{\rule{0.247\textwidth}{0pt}}}
        \subfloat[\label{fig:profit_diff_eps0.025_theta0.01}]{\phantom{\rule{0.247\textwidth}{0pt}}}
    \end{minipage}
    \end{minipage}
    \vspace{-2.em}
    \caption{ 
    Performance comparison for the bilevel problem.}
    \parbox{\textwidth}{\vspace{2mm}
    \footnotesize{\textit{Note.} The evaluations are based on different parameter combinations of risk level ($\epsilon$), radius ($\theta$) parameters, and the number of time steps $T$. Subplots (e)--(p) replicate (a)--(b) but for varying $(\epsilon, \theta)$ combinations. The black horizontal line represents the $3600$s \texttt{TimeLimit}. Dots represent the mean values of the $30$ random runs, with shaded areas representing the $95\%$ percentile interval (from 2.5th to 97.5th). The blue horizontal line indicates a zero difference of the optimal value compared to the proposed SFLA. The \emph{TimeF} (s) of Bonferroni approximation is not displayed because Bonferroni approximation is overly conservative such that it can easily find high \emph{Bilevel Profit} (k\$) (the bilevel objective value) but poor \emph{Actual Profit} (k\$) (plots~\ref{fig:profit_diff_eps0.05_theta0.01} and \ref{fig:profit_diff_eps0.025_theta0.01}). Numbers ending with ``x'' represent the speedup in the average computing time of the proposed SFLA compared to benchmarks. }
    }
    \label{fig:bilevel_results}
\end{figure}

\subsubsection{Results for Bilevel Strategic Bidding}

Figure~\ref{fig:bilevel_results} visualizes the distributions of the evaluation metrics for the proposed SFLA and three benchmarks, i.e., LA, W-CVaR, and Bonferroni. As can be seen from the first two columns of Figure~\ref{fig:bilevel_results}, SFLA demonstrates significant speedup in both \emph{Time} and \emph{TimeF}. In particular, SFLA achieves a speedup of up to $90 \times$
in solving to optimality (\emph{Time}) compared to LA and W-CVaR, and the speedup is up to $40 \times$ in finding the first comparable high-quality solution (\emph{TimeF}) relative to LA and $180\times$ relative to W-CVaR. Although the Bonferroni approximation achieves the fastest computation, its over-conservativeness leads to infeasibility in the majority of parameter settings. It should be noted that the speedup rates compared to LA and W-CVaR decrease for a large $T$, especially for $\epsilon=0.05$; this is because LA and W-CVaR reach the \texttt{TimeLimit} rather than due to a reduction in the computational superiority of the proposed SFLA.

Our Appendix \ref{appen:sup_bilevel} provides supplemental analysis on \emph{N Infeasible} and the actual profits that can be attained by using the proposed SFLA and other benchmarks (the third and fourth columns of Figure~\ref{fig:bilevel_results}). Results demonstrate that the proposed SFLA is more numerically stable and attains comparable actual market profits.



\subsection{Risk Minimization for Robustness Maximization}

\begin{figure}
    \centering    
    \begin{minipage}{\textwidth}
    \vspace{-4mm}
    \includegraphics[width=\columnwidth]{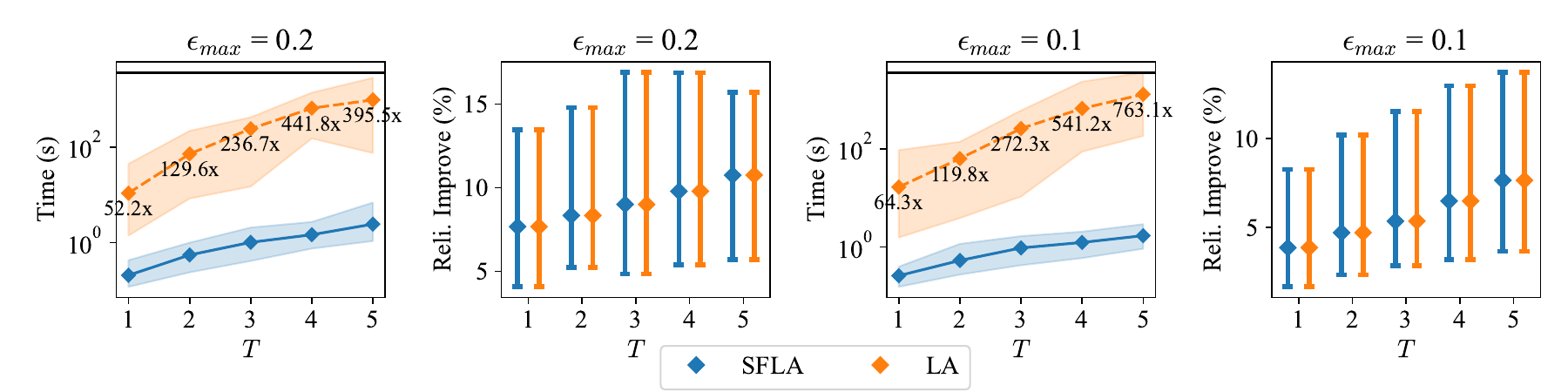}
    \end{minipage}
    \begin{minipage}{\textwidth}
    \hspace{1em}
    \begin{minipage}{\textwidth}
        \centering
        \vspace{-2em}  
        \subfloat[\label{suc_epstheta_time_eps0.2}]{\phantom{\rule{0.247\textwidth}{0pt}}}
        \subfloat[\label{suc_epstheta_time_rel0.2}]{\phantom{\rule{0.247\textwidth}{0pt}}}
        \subfloat[\label{suc_epstheta_time_eps0.1}]{\phantom{\rule{0.247\textwidth}{0pt}}}
        \subfloat[\label{suc_epstheta_time_rel0.1}]{\phantom{\rule{0.247\textwidth}{0pt}}}
    \end{minipage}
    \end{minipage}
    \vspace{-2.5em}
    \caption{ 
    Performance comparison of risk minimization by the proposed SFLA and LA. }
    \parbox{\textwidth}{\vspace{2mm}
    \footnotesize{\textit{Note.} The square dots represent the mean values over 30 random runs, with shaded areas and error bars indicating the 95\% percentile interval. The x-axis represents the length of the optimization horizon $T$. Numbers ending with ``x'' represent the speedup in the average computing time of the proposed SFLA compared to LA. The first two plots correspond to results where the upper bound of the risk level $\epsilon$ is set to 0.2, while the remaining two use $\epsilon = 0.1$. In the original UC problem, the Wasserstein radius $\theta$ is set to 0.1, and the risk levels are set to the upper bounds mentioned above. \emph{Cost Tol (\%)} is set to 3\%. }
    \vspace{3mm}}
    \label{fig:suc_epstheta}
\end{figure}

This section compares the proposed SFLA and LA extended for risk minimization. ExactS, WCVaR, and Bonferroni are not compared as their limitations have been illustrated in previous sections. The test case considered is the UC problem in Section~\ref{sec:UC}. We first solve the chance-constrained UC with SFLA at $\epsilon_{\max}$ to obtain the base cost, and then set the bearable cost $\tau$ in Problem~\eqref{eq:eps_robust_satisfi} to the base cost multiplied by $1+\textit{Cost Tol (\%)}$. Regarding the computing time, as can be seen in Figure~\ref{suc_epstheta_time_eps0.2}, the proposed SFLA leads to significant computational speedup over LA, and the speedup increases from $52.2\times$ to more than $400\times$ for $\epsilon_{\max} =0.2$ as the problem scales up (a greater optimization horizon $T$). Moreover, when $\epsilon_{\max} =0.1$ in Figure~\ref{suc_epstheta_time_eps0.1}, the speedup is up to $763\times$. This is because a smaller $\epsilon_{\max} =0.1$ leads to a stronger strengthening. Finally, as can be seen in Figures \ref{suc_epstheta_time_rel0.2} and \ref{suc_epstheta_time_rel0.1}, both SFLA and LA improve the problem robustness (\emph{Reli. Improve} (\%), which represents the absolute improvement of the out-of-sample \emph{Reli.} (\%) after robustness maximization). The difference in the mean values of the \emph{Reli. Improve (\%)} is less than 0.01\% (the negligible difference is attributed to the existence of multiple optimal solutions). This supports our theoretical results that the extended SFLA is still valid and brings computational speedup compared to LA without being more conservative, even when applied to risk minimization. Our Appendix~\ref{res:risk_vs_radius} further shows the radius maximization results for robustness maximization.

\section{Conclusion}
\label{sec:conclusion}

This work proposes a novel Strengthened and Faster Linear Approximation (SFLA) for Wasserstein distributionally robust joint chance constraints under right-hand-side uncertainty (RHS-WDRJCC). The proposed SFLA enables more efficient and less conservative solutions to a class of real-world large-scale chance-constrained optimization problems, such as chance-constrained unit commitment (UC) and bilevel optimization for strategic bidding in electricity markets with chance-constrained economic dispatch. By applying valid inequalities to an existing linear approximation (LA), the proposed SFLA reduces the number of constraints and tightens the feasible region for the ancillary variables, resulting in a significant computational speedup. Despite this tightening, we demonstrate that the proposed SFLA can even be a less conservative inner approximation compared to LA and the worst-case conditional value-at-risk (W-CVaR) method. We also extend the proposed SFLA to robustness maximizing, \textcolor{black}{a decision-making paradigm that can be more interpretable to set the risk level and the Wasserstein radius}. We prove the equivalence of the feasible regions between risk minimization and radius maximization as two formulations for robustness maximization, and show that radius maximization can be both more computationally efficient and more robust.

Extensive numerical experiments demonstrate the superiority and significance of the proposed SFLA. In the UC problem, we demonstrate that the proposed SFLA achieves up to $10\times$ and on average $3.8\times$ computational speedup compared to the existing exact and strengthened MIP reformulation (ExactS) in finding comparable high-quality solutions. We also demonstrate that the proposed SFLA achieves approximation quality comparable to the exact reformulation ExactS, except for minor inferiority only under the combination of a high risk $\epsilon$, a large number of historical data ($N$), and a small ambiguity set radius $\theta$. In the bilevel strategic bidding problem where the exact formulation is not applicable due to its non-convexity, the proposed SFLA can lead to $90 \times$ computational speedup than existing linear approximation methods (LA and W-CVaR). In robustness maximization, the proposed SFLA demonstrated over $100\times$ computational speedup compared to other convex approximations with the same solution quality.

It is worth noting that the applications of the proposed SFLA are not limited to the UC and bilevel energy bidding problems demonstrated in this paper. The method is also applicable to other large-scale, real-world problems that are already complex in their deterministic form and require tractable approaches to handle uncertainty. Examples include large-scale transportation problems, multilevel market-clearing problems, and multilevel infrastructure planning that involves a lower-level market being cleared before uncertainty is realized. 

Historical data plays a significant role in affecting the complexity and solution quality. Future work could also explore scenario reduction techniques~\citep{arpon2018scenario}
to further improve SFLA.


\ACKNOWLEDGMENT{
This paper was supported by the Engineering and Physical Sciences Research Council (EPSRC) (grant reference number EP/W027321/1). Yihong Zhou's and Yuxin Xia's works were also supported by the Engineering Studentship from the University of Edinburgh. 
Hanbin Yang's work is funded by the National Natural Science Foundation of China under Grant 72201232 and 72231008, Guangdong Provincial Key Laboratory of Mathematical Foundations for Artificial Intelligence (2023B1212010001), and Shenzhen Key Laboratory
of Crowd Intelligence Empowered Low-Carbon Energy Network, China (project number ZDSYS20220606100601002).
%
}

%
%
%


\bibliographystyle{informs2014} 
\bibliography{bibref} 

@article{luedtke2010integer,
  title={An integer programming approach for linear programs with probabilistic constraints},
  author={Luedtke, James and Ahmed, Shabbir and Nemhauser, George L},
  journal={Mathematical programming},
  volume={122},
  number={2},
  pages={247--272},
  year={2010},
  publisher={Springer}
}

@article{arpon2018scenario,
  title={Scenario reduction for stochastic programs with conditional value-at-risk},
  author={Arp{\'o}n, Sebasti{\'a}n and Homem-de-Mello, Tito and Pagnoncelli, Bernardo},
  journal={Mathematical Programming},
  volume={170},
  number={1},
  pages={327--356},
  year={2018},
  publisher={Springer}
}

@misc{gurobi,
  author = {{Gurobi Optimization, LLC}},
  title = {{Gurobi Optimizer Reference Manual}},
  year = 2024,
  url = "https://www.gurobi.com"
}

@article{CCP_power_system_review,
  title={Data-driven decision making in power systems with probabilistic guarantees: Theory and applications of chance-constrained optimization},
  author={Geng, Xinbo and Xie, Le},
  journal={Annual reviews in control},
  volume={47},
  year={2019},
  publisher={Elsevier}
}

@book{Stochastic-programming,
  title={Stochastic programming},
  author={Pr{\'e}kopa, Andr{\'a}s},
  volume={324},
  year={2013},
  publisher={Springer Science \& Business Media}
}

@ARTICLE{yang2019analytical,  author={Yang, Yue and Wu, Wenchuan and Wang, Bin and Li, Mingjie},  journal={IEEE Transactions on Power Systems},   title={Analytical Reformulation for Stochastic Unit Commitment Considering Wind Power Uncertainty With Gaussian Mixture Model},   year={2020},  volume={35},  number={4},  pages={2769-2782}}

@article{exact_milp_strengthened,
  title={Distributionally robust chance-constrained programs with right-hand side uncertainty under {W}asserstein ambiguity},
  author={Ho-Nguyen, Nam and K{\i}l{\i}n{\c{c}}-Karzan, Fatma and K{\"u}{\c{c}}{\"u}kyavuz, Simge and Lee, Dabeen},
  journal={Mathematical Programming},
  year={2021},
  publisher={Springer}
}

@ARTICLE{7828161,
  author={Wang, Yishen and Dvorkin, Yury and Fernández-Blanco, Ricardo and Xu, Bolun and Qiu, Ting and Kirschen, Daniel S.},
  journal={IEEE Transactions on Sustainable Energy}, 
  title={Look-Ahead Bidding Strategy for Energy Storage}, 
  year={2017},
  volume={8},
  number={3},
  pages={1106-1117},}

@article{2022_OR_exact_DRO,
  title={Data-driven chance constrained programs over {W}asserstein balls},
  author={Chen, Zhi and Kuhn, Daniel and Wiesemann, Wolfram},
  journal={Operations Research},
  year={2022},
  publisher={INFORMS}
}

@article{chen2023approximations,
  title={On approximations of data-driven chance constrained programs over {W}asserstein balls},
  author={Chen, Zhi and Kuhn, Daniel and Wiesemann, Wolfram},
  journal={Operations Research Letters},
  volume={51},
  number={3},
  pages={226--233},
  year={2023},
  publisher={Elsevier}
}

@article{TOMASSON2020114251,
title = {Optimal offer-bid strategy of an energy storage portfolio: A linear quasi-relaxation approach},
journal = {Applied Energy},
volume = {260},
pages = {114251},
year = {2020},
issn = {0306-2619},
author = {Egill Tómasson and Mohammad Reza Hesamzadeh and Frank A. Wolak}
}

@article{DIMITRIADIS2022123026,
title = {Strategic bidding of an energy storage agent in a joint energy and reserve market under stochastic generation},
journal = {Energy},
volume = {242},
pages = {123026},
year = {2022},
issn = {0360-5442},
author = {Christos N. Dimitriadis and Evangelos G. Tsimopoulos and Michael C. Georgiadis},
}

@ARTICLE{8036231,
  author={Nasrolahpour, Ehsan and Kazempour, Jalal and Zareipour, Hamidreza and Rosehart, William D.},
  journal={IEEE Transactions on Sustainable Energy}, 
  title={A Bilevel Model for Participation of a Storage System in Energy and Reserve Markets}, 
  year={2018},
  volume={9},
  number={2},
  pages={582-598}}

@article{ordoudis2021energy,
  title={Energy and reserve dispatch with distributionally robust joint chance constraints},
  author={Ordoudis, Christos and Nguyen, Viet Anh and Kuhn, Daniel and Pinson, Pierre},
  journal={Operations Research Letters},
  volume={49},
  number={3},
  pages={291--299},
  year={2021},
  publisher={Elsevier}
}

@book{scarf1957min,
  title={A min-max solution of an inventory problem},
  author={Scarf, Herbert E and Arrow, KJ and Karlin, S},
  year={1957},
  publisher={Rand Corporation Santa Monica}
}

@article{delage2010distributionally,
  title={Distributionally robust optimization under moment uncertainty with application to data-driven problems},
  author={Delage, Erick and Ye, Yinyu},
  journal={Operations research},
  volume={58},
  number={3},
  pages={595--612},
  year={2010},
  publisher={INFORMS}
}

@article{gao2023distributionally,
  title={Distributionally robust stochastic optimization with {W}asserstein distance},
  author={Gao, Rui and Kleywegt, Anton},
  journal={Mathematics of Operations Research},
  volume={48},
  number={2},
  pages={603--655},
  year={2023},
  publisher={INFORMS}
}

@article{mohajerin2018data,
  title={Data-driven distributionally robust optimization using the {W}asserstein metric: Performance guarantees and tractable reformulations},
  author={Mohajerin Esfahani, Peyman and Kuhn, Daniel},
  journal={Mathematical Programming},
  volume={171},
  number={1},
  pages={115--166},
  year={2018},
  publisher={Springer}
}

@article{ding2022distributionally,
  title={Distributionally robust joint chance-constrained optimization for networked microgrids considering contingencies and renewable uncertainty},
  author={Ding, Yifu and Morstyn, Thomas and McCulloch, Malcolm D},
  journal={IEEE Transactions on Smart Grid},
  year={2022},
  publisher={IEEE}
}

@article{van2018exact,
  title={An exact solution method for the hydrothermal unit commitment under wind power uncertainty with joint probability constraints},
  author={van Ackooij, Wim and Finardi, Erlon Cristian and Ramalho, Guilherme Matiussi},
  journal={IEEE Transactions on Power Systems},
  volume={33},
  number={6},
  pages={6487--6500},
  year={2018},
  publisher={IEEE}
}

@article{jiang2022also,
  title={{ALSO-X} and {ALSO-X}+: Better convex approximations for chance constrained programs},
  author={Jiang, Nan and Xie, Weijun},
  journal={Operations Research},
  volume={70},
  number={6},
  pages={3581--3600},
  year={2022},
  publisher={INFORMS}
}

@ARTICLE{MISO_UC_difficulty,
  author={Chen, Yonghong and Casto, Aaron and Wang, Fengyu and Wang, Qianfan and Wang, Xing and Wan, Jie},
  journal={IEEE Transactions on Power Systems}, 
  title={Improving Large Scale Day-Ahead Security Constrained Unit Commitment Performance}, 
  year={2016},
  volume={31},
  number={6},
  pages={4732-4743}}

@article{sun2018novel,
  title={A novel decomposition and coordination approach for large day-ahead unit commitment with combined cycle units},
  author={Sun, Xiaorong and Luh, Peter B and Bragin, Mikhail A and Chen, Yonghong and Wan, Jie and Wang, Fengyu},
  journal={IEEE Transactions on Power Systems},
  volume={33},
  number={5},
  pages={5297--5308},
  year={2018},
  publisher={IEEE}
}

@article{beck2023survey,
  title={A survey on bilevel optimization under uncertainty},
  author={Beck, Yasmine and Ljubi{\'c}, Ivana and Schmidt, Martin},
  journal={European Journal of Operational Research},
  volume={311},
  number={2},
  pages={401--426},
  year={2023},
  publisher={Elsevier}
}

@article{wu2016solution,
  title={A solution to the chance-constrained two-stage stochastic program for unit commitment with wind energy integration},
  author={Wu, Zhi and Zeng, Pingliang and Zhang, Xiao-Ping and Zhou, Qinyong},
  journal={IEEE Transactions on Power Systems},
  volume={31},
  number={6},
  pages={4185--4196},
  year={2016},
  publisher={IEEE}
}

@article{wang2016risk,
  title={Risk adjustable day-ahead unit commitment with wind power based on chance constrained goal programming},
  author={Wang, Yang and Zhao, Shuqiang and Zhou, Zhi and Botterud, Audun and Xu, Yan and Chen, Runze},
  journal={IEEE Transactions on Sustainable Energy},
  volume={8},
  number={2},
  pages={530--541},
  year={2016},
  publisher={IEEE}
}

@article{chen2007robust,
  title={A robust optimization perspective on stochastic programming},
  author={Chen, Xin and Sim, Melvyn and Sun, Peng},
  journal={Operations research},
  volume={55},
  number={6},
  pages={1058--1071},
  year={2007},
  publisher={INFORMS}
}

@misc{xu2017application,
  title={Application of large-scale synthetic power system models for energy economic studies},
  author={Xu, Ti and Birchfield, Adam B and Gegner, Kathleen M and Shetye, Komal S and Overbye, Thomas J},
  year={2017},
  journal={ }
}

@article{heitsch2022convex,
  title={On convex lower-level black-box constraints in bilevel optimization with an application to gas market models with chance constraints},
  author={Heitsch, Holger and Henrion, Ren{\'e} and Kleinert, Thomas and Schmidt, Martin},
  journal={Journal of Global Optimization},
  volume={84},
  number={3},
  pages={651--685},
  year={2022},
  publisher={Springer}
}

@article{paredes2023stacking,
  title={Stacking revenues from flexible DERs in multi-scale markets using tri-level optimization},
  author={Paredes, {\'A}ngel and Aguado, Jos{\'e} A and Essayeh, Chaimaa and Xia, Yuxin and Savelli, Iacopo and Morstyn, Thomas},
  journal={IEEE Transactions on Power Systems},
  volume={39},
  number={2},
  pages={3949--3961},
  year={2023},
  publisher={IEEE}
}

@article{schwartz2011makes,
  title={What makes a good decision? Robust satisficing as a normative standard of rational decision making},
  author={Schwartz, Barry and Ben-Haim, YAKOV and Dacso, Cliff},
  journal={Journal for the Theory of Social Behaviour},
  volume={41},
  number={2},
  pages={209--227},
  year={2011},
  publisher={Wiley Online Library}
}

@article{zhang2025integer,
  title={Integer programming approaches for distributionally robust chance constraints with adjustable risks},
  author={Zhang, Yiling},
  journal={Computers \& Operations Research},
  volume={177},
  pages={106974},
  year={2025},
  publisher={Elsevier}
}

@article{long2023robust,
  title={Robust satisficing},
  author={Long, Daniel Zhuoyu and Sim, Melvyn and Zhou, Minglong},
  journal={Operations Research},
  volume={71},
  number={1},
  pages={61--82},
  year={2023},
  publisher={INFORMS}
}

@article{jiang2024also,
  title={{ALSO-X}\#: Better convex approximations for distributionally robust chance constrained programs},
  author={Jiang, Nan and Xie, Weijun},
  journal={Mathematical Programming},
  pages={1--64},
  year={2024},
  publisher={Springer}
}

@article{zhao2023distributionally,
  title={Distributionally robust chance-constrained p-hub center problem},
  author={Zhao, Yue and Chen, Zhi and Zhang, Zhenzhen},
  journal={INFORMS Journal on Computing},
  volume={35},
  number={6},
  pages={1361--1382},
  year={2023},
  publisher={INFORMS}
}

@article{gabrel2014recent,
  title={Recent advances in robust optimization: An overview},
  author={Gabrel, Virginie and Murat, C{\'e}cile and Thiele, Aur{\'e}lie},
  journal={European journal of operational research},
  volume={235},
  number={3},
  pages={471--483},
  year={2014},
  publisher={Elsevier}
}

@misc{zhao2024viabilitystochasticeconomicdispatch,
      title={On the Viability of Stochastic Economic Dispatch for Real-Time Energy Market Clearing}, 
      author={Haoruo Zhao and Mathieu Tanneau and Pascal Van Hentenryck},
      year={2024},
      eprint={2308.06386},
      archivePrefix={arXiv},
      primaryClass={math.OC},
      url={https://arxiv.org/abs/2308.06386}, 
}

@article{jiang2024terminator,
  title={The Terminator: An Integration of Inner and Outer Approximations for Solving {W}asserstein Distributionally Robust Chance Constrained Programs via Variable Fixing},
  author={Jiang, Nan and Xie, Weijun},
  journal={INFORMS Journal on Computing},
  year={2024},
  publisher={INFORMS}
}

@misc{Zhou2024_code,
  author =     {Zhou, Yihong and Xia, Yuxin and Yang, Hanbin and Morstyn, Thomas},
  publisher =  {INFORMS Journal on Computing},
  title =      {Strengthened and Faster Linear Approximation to Joint Chance Constraints with Wasserstein Ambiguity},
  year =       {2026},
  doi =        {10.1287/ijoc.2024.1073.cd},
  note =       {Available for download at https://github.com/INFORMSJoC/2024.1073},
}

\newpage
\renewcommand{\theHsection}{A\arabic{section}} 
\begin{APPENDIX}{Supplemental Information}\label{appendix}

\section{Supplemental Formulations}
\label{appen:exactS}

The MIP reformulation of the exact reformulation \eqref{eq:wdrjcc_abstract} is provided below: 
\begin{subequations}\label{eq:MILP_constr1}
    \begin{align}
        & \b z \in \{0,1\}^N, s \geq 0, \b r \geq \b 0, \label{con:s_r_z}\\
        & \epsilon N s - \sum\limits_{i=1}^N r_i \geq \theta N,  \label{con:proportion} \\
        & M(1-z_i) \geq s-r_i, &&\hspace{-50pt} \forall i\in [N], \label{con:bigM} \\
        & \dfrac{\b b_p^\top \b \xi_i + d_p - \b a_p^\top \b x}{\| \b b_p \|_*} + M z_i \geq s-r_i, &&\hspace{-50pt} \forall i\in [N], p\in[P].
    \end{align}
\end{subequations}
The strengthened exact MIP reformulation (ExactS), presented as Formulation (20) in~\citep{exact_milp_strengthened}, is reproduced below for convenience:
\begin{subequations}\label{eq:ExactS}
    \begin{align}
        & \b z, s, \b r \text{ satisfy \eqref{con:s_r_z}, \eqref{con:proportion}, \eqref{con:bigM} and} \sum_{i\in[N]} z_i \leq \lfloor \epsilon N \rfloor, \label{eq:exacts1}  \\
        & \dfrac{\b b_p^\top \b \xi_i + d_p - \b a_p^\top \b x}{\| \b b_p \|_*} +  \dfrac{-\b b_p^\top \b \xi_i + q_p}{\| \b b_p \|_*} z_i \geq s - r_i, &&\hspace{-50pt} \forall i \in [N]_p, p\in [P], \\
        & \dfrac{q_p + d_p - \b a_p^\top \b x}{\| \b b_p \|_*} \geq s, &&\hspace{-50pt} \forall p\in [P].\label{eq:exacts2}
    \end{align}
\end{subequations}
where all the notations have been defined in the main content. The $\b x$-feasible region of the ExactS can be expressed as $\mathcal{X}_\text{ExactS} = \left\{ \b x \in \mathcal{X} \mid  \exists\ \b z, s, \b r :~\eqref{eq:exacts1}\mbox{--}\eqref{eq:exacts2} \right\}$.

{\color{black}
\section{Proofs}
\label{appen:proof}

Completing our proofs requires the introduction of the following lemmas:
\begin{lemma}[Lemma 1 in \cite{exact_milp_strengthened}]
\label{lemma:existence}
    For any $\b x \in \mathcal{X}_\text{Exact}$, there exists ($\b r$, $s$) such that $s$ is equal to the (k+1)-th smallest element in $\{ \mathrm{dist}(\b \xi_i, \mathcal{S}(\b x)) \}_{i\in [N]}$ and constraints~\eqref{eq:s_r}--\eqref{eq:dist_constr} are satisfied.
\end{lemma}

\begin{lemma}\label{lemma:cvar}
    $\mathcal{X}_\text{LA}(\b\kappa)$ defined through~\eqref{eq:LP_constr_1} is equivalent to the following set:
    \begin{align}\label{eq:feasi_set_inner_cvar}
        \left\{ \b x \in \mathcal{X} \;\mid\; \frac{\theta}{\epsilon} + \mathbb{P}_N\mhyphen\mathrm{CVaR_{1-\epsilon}}(-\widehat{\mathrm{dist}}(\b \xi, \mathcal{S}(\b x))) \leq 0 \right\},
    \end{align}
    where $\mathbb{P}_N\mhyphen\mathrm{CVaR_{1-\epsilon}}(v(\b \xi)) \coloneqq \min_{s'} \left\{ s' + \frac{1}{\epsilon N} \sum_{i \in [N]} \max \{0, v(\b \xi_i) - s'\} \right\}$ is the CVaR of the random variable $v(\b \xi)$ under the empirical probability distribution $\mathbb{P}_N$~\citep{Stochastic-programming}.
\end{lemma}
\proof{proof}
   \eqref{eq:feasi_set_inner_cvar} $\Longrightarrow$~\eqref{eq:LP_constr_1}: Denote $-\widehat{\mathrm{dist}}(\b \xi_i, \mathcal{S}(\b x))$ by $v(\b \xi_i)$, assume $s'$ is selected such that $s' + \frac{1}{\epsilon N} \sum\limits_{i \in [N]} \max \{0, v(\b \xi_i) - s'\}$ is minimized, and take $r_i^* = \max \{0, v(\b \xi_i) - s'\}$ (so $r_i^* \geq 0$). The inequality in~\eqref{eq:feasi_set_inner_cvar} indicates:
    \begin{equation}
         \dfrac{\theta}{\epsilon} + s' + \frac{1}{\epsilon N} \sum_{i=1}^N r_i^* \leq 0 
         \Leftrightarrow 
         \epsilon N s' + \sum_{i=1}^N r_i^* \leq -\theta N. \label{eq:right_direction}
    \end{equation}    
    Taking $s^* = -s'$, now we will prove that the tuple $(\b x, s^*, \b r^*)$ satisfy all the constraints in~\eqref{eq:LP_constr_1}, namely~\eqref{eq:LP1_s_r}--\eqref{eq:LP1_last}, for any $\b x$ belonging to the set~\eqref{eq:feasi_set_inner_cvar}. Since $r_i^* \geq 0$ and $-\theta N \leq 0$,~\eqref{eq:right_direction} implies $s' \leq 0$ (i.e., $s^* \geq 0$) and the following inequality $\epsilon N s - \sum_{i=1}^N r_i \geq \theta N$,
    which means~\eqref{eq:LP1_s_r} and~\eqref{eq:LP1_epsilonNs} are satisfied. It remains to prove the satisfaction of~\eqref{eq:LP1_last}. Based on the setting of $r_i^*$, we have:
    \begin{equation}
    \begin{split}
        s^* - r_i^* & = s^* - \max \left\{0, v(\b \xi_i) + s^* \right\} = s^* - \max \left\{0, s^* -\widehat{\mathrm{dist}}(\b \xi_i, \mathcal{S}(\b x)) \right\}.
    \end{split}
    \label{eq:s-ri}
    \end{equation}
    Notice that constraint~\eqref{eq:LP1_last} is equivalent to:
    \begin{equation}
       \widehat{\mathrm{dist}}(\b \xi_i, \mathcal{S}(\b x)) = \min_{p \in [P]} \left\{\kappa_i \frac{\b b_p^\top \b \xi + d_p - \b a_p^\top \b x}{\| \b b_p \|_*} \right\}\geq s-r_i \quad \forall i\in [N]. \label{eq:LP1_last_eq}
    \end{equation}
    When $\widehat{\mathrm{dist}}(\b \xi_i, \mathcal{S}(\b x)) \geq s^*$,~\eqref{eq:LP1_last} must be satisfied for $(\b x, s^*, \b r^*)$ due to $r_i^* \geq 0$. When $\widehat{\mathrm{dist}}(\b \xi_i, \mathcal{S}(\b x)) < s^*$, we have $s^*-r_i^* = \widehat{\mathrm{dist}}(\b \xi_i, \mathcal{S}(\b x))$ based on~\eqref{eq:s-ri}. Therefore,~\eqref{eq:LP1_last_eq} and equivalently~\eqref{eq:LP1_last} is satisfied.

   \eqref{eq:feasi_set_inner_cvar} $\Longleftarrow$~\eqref{eq:LP_constr_1}: For any $\b x$ belonging to the $\mathcal{X}_\text{LA}(\b\kappa)$, by definition there exist $(\b r, s)$ such that constraints~\eqref{eq:LP1_s_r}--\eqref{eq:LP1_last}, are satisfied. Among these constraints,~\eqref{eq:LP1_last} (equivalently~\eqref{eq:LP1_last_eq}) informs that $r_i \geq s - \widehat{\mathrm{dist}}(\b \xi_i, \mathcal{S}(\b x))$.
    Combining with $r_i \geq 0$ in~\eqref{eq:LP1_s_r}, we further have $r_i \geq \max \left\{0, s - \widehat{\mathrm{dist}}(\b \xi_i, \mathcal{S}(\b x)) \right\}$.
    Then we have $\epsilon N s - \sum_{i \in [N]} r_i \leq \epsilon N s  - \sum_{i \in [N]} \max \left\{0, s - \widehat{\mathrm{dist}}(\b \xi_i, \mathcal{S}(\b x)) \right\}$.
    Because of $\epsilon N s - \sum\limits_{i\in [N]}r_i \geq \theta N$ in~\eqref{eq:LP1_epsilonNs}, we have $\epsilon N s  - \sum_{i \in [N]} \max \left\{0, s - \widehat{\mathrm{dist}}(\b \xi_i, \mathcal{S}(\b x)) \right\} \geq \theta N$,
    which means $\frac{\theta}{\epsilon} - s + \frac{1}{\epsilon N}\sum_{i \in [N]} \max \left\{0, s - \widehat{\mathrm{dist}}(\b \xi_i, \mathcal{S}(\b x)) \right\} \leq 0$.
    Now take $s' = -s$, we further have $\frac{\theta}{\epsilon} + s' + \frac{1}{\epsilon N}\sum_{i \in [N]} \max \left\{0, -s' - \widehat{\mathrm{dist}}(\b \xi_i, \mathcal{S}(\b x)) \right\} \leq 0$, 
    which leads to $\frac{\theta}{\epsilon} + \min_{s'}\left\{ s' + \frac{1}{\epsilon N}\sum_{i \in [N]} \max\{0, -s' + v(\b \xi_i) \right\} \leq 0$.
    Therefore, the constraint in~\eqref{eq:feasi_set_inner_cvar} is satisfied.
$\hfill\blacksquare$\endproof

\begin{lemma}\label{lemma:our_exist}
    For any $\b x \in \mathcal{X}_\text{LA}(\b\kappa)$, there exists ($\b r$, $s$) such that $s$ is equal to the $(k+1)$-th smallest value amongst $\{ \widehat{\mathrm{dist}}(\b \xi_i, \mathcal{S}(\b x))\}_{i\in [N]}$ and the constraints that define $\mathcal{X}_\text{LA}(\b \kappa)$, namely~\eqref{eq:LP1_s_r}--\eqref{eq:LP1_last}, are satisfied.
\end{lemma}
\proof{proof}
    Notice that the CVaR of a random variable $v(\b \xi)$ with $\b \xi \sim \mathbb{P}_N$ (defined in Lemma \ref{lemma:cvar}) can be equivalently expressed as the following primal and dual forms:  
    \begin{subequations}\label{eq:CVaR_optional}
    \begin{align}
    \mbox{Primal:\quad}
    \mathbb{P}_N\mhyphen\mathrm{CVaR_{1-\epsilon}}(v(\b \xi)) = 
        \min_{s'\in \mathbb{R}, \b r \geq 0} \quad &s' + \frac{1}{\epsilon N} \sum\nolimits_{i=1}^N r_i \label{eq:CVaR_optional_obj}\\
        \mbox{s.t.} \quad & r_i \geq v(\b \xi_i ) - s',\quad i\in [N].
    \end{align}
    \end{subequations}
    \vspace{-10mm}
    \begin{subequations}
    \begin{align}
    \mbox{Dual:\quad}
        \mathbb{P}_N\mhyphen\mathrm{CVaR_{1-\epsilon}}(v(\b \xi)) = 
        \max_{\b 0\leq \b y \leq \b 1} \quad & \frac{1}{\epsilon N} \sum\nolimits_{i\in [N]} v(\b \xi_i) y_i\\
        \mbox{s.t.} \quad &\frac{1}{N} \sum\nolimits_{i\in [N]} y_i = \epsilon.
    \end{align} 
    \end{subequations}     
    Without loss of generality, assume that we have an ordering $v(\b \xi_1) \geq \cdots \geq v(\b \xi_N)$. Then a primal-dual optimal pair for CVaR is given by $s'^* = v(\b \xi_{k+1})$, $r_i^* = \max\{ 0, v(\b \xi_i) - s'^* \}$, and 
    $y_i^* = 
    \begin{cases} 
      1, & i = 1,\cdots , k \\
      \epsilon N - k, & i=k+1 \\
      0, & i > k+1.
   \end{cases} $.    
    The validity of the above primal-dual optimal pair holds because of the satisfaction of primal and dual constraints ($\b y^* \in [0,1]$ and $\b r^* \geq 0$) and the equality of the objectives of the primal and dual problems $\frac{1}{\epsilon N} \sum_{i\in [N]} v(\b \xi_i) y_i^*
        = s'^* + \frac{1}{\epsilon N} \sum_{i \in [N]} r_i^*$,
    which holds because of $\frac{1}{\epsilon N} \sum_{i\in [N]} v(\b \xi_i) y_i^*
         =
        \frac{1}{\epsilon N} \left( \sum_{i\in [k]} v(\b \xi_i) + v(\b \xi_{k+1}) (\epsilon N - k) \right)
         = v(\b \xi_{k+1}) + \frac{1}{\epsilon N} \left( \sum_{i\in [k]} v(\b \xi_i) - k v(\b \xi_{k+1}) \right)$
    and
    $
        s'^* + \frac{1}{\epsilon N} \sum_{i \in [N]} r_i^*  = v(\b \xi_{k+1}) + \frac{1}{\epsilon N} \sum_{i \in [N]} \max\{ 0, v(\b \xi_i) - v(\b \xi_{k+1}) \} 
        = v(\b \xi_{k+1}) + \frac{1}{\epsilon N} \sum_{i \in [k]} \left( v(\b \xi_i) - v(\b \xi_{k+1}) \right) 
        = v(\b \xi_{k+1}) + \frac{1}{\epsilon N} \left( \sum_{i\in [k]} v(\b \xi_i) - k v(\b \xi_{k+1}) \right)
    $, where the reduction from $\sum_{i \in [N]}$ to $\sum_{i \in [k]}$ is valid because of our ordering of $v(\b \xi_i)$ such that $v(\b \xi_i)-v(\b \xi_{k+1})\leq 0$ for $i\geq k+1$.
    
    Now take $v(\b \xi_i) = - \widehat{\mathrm{dist}}(\b \xi_i, \mathcal{S}(\b x))$. Based on the CVaR interpretation of $\mathcal{X}_\text{LA}(\b\kappa)$ in~\eqref{eq:feasi_set_inner_cvar} and the CVaR value expressed as~\eqref{eq:CVaR_optional_obj}, for any $\b x \in \mathcal{X}_\text{LA}(\b\kappa)$, we have:    
    \begin{align*}
        \frac{\theta}{\epsilon} + \mathbb{P}_N\mhyphen\mathrm{CVaR_{1-\epsilon}}\left(-\widehat{\mathrm{dist}}(\b \xi_i, \mathcal{S}(\b x))\right) \leq 0 
        \Leftrightarrow 
        \frac{\theta}{\epsilon} + s'^* + \frac{1}{\epsilon N} \sum_{i=1}^N r_i^* \leq 0 
        \Leftrightarrow  
        \epsilon N s'^* + \sum_{i=1}^N r_i^* \leq -\theta N.
    \end{align*}      
    Since $r_i^* \geq 0$ and $-\theta N \leq 0$, the above result indicates $s^* \coloneqq -s'^* \geq 0$ and $\epsilon N s^* - \sum_{i \in [N]} r_i^* \geq \theta N$,
    which means~\eqref{eq:LP1_s_r} and~\eqref{eq:LP1_epsilonNs} are satisfied. It remains to prove the satisfaction of~\eqref{eq:LP1_last}. Based on the setting of $r_i^*$, we have $s^* - r_i^* = s^* - \max \{0, v(\b \xi_i) + s^*\} = s^* - \max \{0, s^* -\widehat{\mathrm{dist}}(\b \xi_i, \mathcal{S}(\b x))\}$.
    Recall (see~\eqref{eq:LP1_last_eq}) that constraint~\eqref{eq:LP1_last} is equivalent to $\widehat{\mathrm{dist}}(\b \xi_i, \mathcal{S}(\b x)) = \min_{p \in [P]} \left\{\kappa_i \frac{\b b_p^\top \b \xi + d_p - \b a_p^\top \b x}{\| \b b_p \|_*} \right\}\geq s-r_i \quad i\in [N]$.
    When $\widehat{\mathrm{dist}}(\b \xi_i, \mathcal{S}(\b x)) \geq s^*$,~\eqref{eq:LP1_last} must be satisfied due to $r_i^* \geq 0$. When $\widehat{\mathrm{dist}}(\b \xi_i, \mathcal{S}(\b x)) < s^*$, we have $s^*-r_i^* = \widehat{\mathrm{dist}}(\b \xi_i, \mathcal{S}(\b x))$, so~\eqref{eq:LP1_last} is satisfied as well. Therefore, we have proved that the pair $(\b x, s^*, \b r^*)$ satisfies constraints in~\eqref{eq:LP_constr_1} for any $\b x \in \mathcal{X}_\text{LA}(\b\kappa)$, with $s^* = -s'^* = -v(\b \xi_{k+1})$ and $r_i^* = \max\{ 0, v(\b \xi_i) + s^* \}$. Finally, based on our ordering, $v(\b \xi_{k+1}) =-\widehat{\mathrm{dist}}(\b \xi_{k+1}, \mathcal{S}(\b x))$ is the $(k+1)$-th greatest value amongst $\{ v(\b \xi_i) \}_{i \in [N]}$, so $s^* = -v(\b \xi_{k+1})$ is the $(k+1)$-th smallest value amongst $\{\widehat{\mathrm{dist}}(\b \xi_i, \mathcal{S}(\b x))\}_{i\in [N]}$.
$\hfill\blacksquare$\endproof

\subsection*{Proof of Proposition~\ref{proposition}}

\proof{proof}
    It suffices to show that replacing the constraint~\eqref{eq:dist_constr} with constraints~\eqref{eq:MILP_reduced_part}--\eqref{eq:MILP_reduced_q} has no impact on the $\b x$-feasible region.
    
    \noindent \eqref{eq:MILP_constr_reduced}$\Longrightarrow$\eqref{eq:wdrjcc_abstract}: Based on the definition of $[N]_p$, we have $\b b_p^\top \b \xi_i \geq q_p, \  \forall i \in [N] \setminus [N]_p,\ p\in [P]$.
    Then, based on~\eqref{eq:MILP_reduced_q}, for any $i \in [N] \setminus [N]_p$ with $p\in [P]$, we have:
    $
    \Big(\dfrac{\b b_p^\top \b \xi_i + d_p - \b a_p^\top \b x}{\| \b b_p \|_*}\Big)^+ \geq \dfrac{\b b_p^\top \b \xi_i + d_p - \b a_p^\top \b x}{\| \b b_p \|_*} \geq  \dfrac{q_p + d_p - \b a_p^\top \b x}{\| \b b_p \|_*} \geq s
    $.
    As $r_i \geq 0$, we further have: $\Big(\dfrac{\b b_p^\top \b \xi_i + d_p - \b a_p^\top \b x}{\| \b b_p \|_*}\Big)^+ \geq s \geq s - r_i,$
    which, combined with~\eqref{eq:MILP_reduced_part}, implies the satisfaction of~\eqref{eq:dist_constr}.

    \noindent \eqref{eq:MILP_constr_reduced}$\Longleftarrow$\eqref{eq:wdrjcc_abstract}: Based on Lemma \ref{lemma:existence}, for any $\b x \in \mathcal{X}_\text{Exact}$ defined though~\eqref{eq:wdrjcc_abstract}, there always exists ($\b r^*$, $s^*$) such that $s^*$ is equal to the (k+1)-th smallest value amongst the set $\{ \mathrm{dist}(\b \xi_i, \mathcal{S}(\b x)) \}_{i\in [N]}$ and constraints~\eqref{eq:s_r}--\eqref{eq:dist_constr} are satisfied. Without loss of generality, we assume the following ordering: $\mathrm{dist}(\b \xi_1, \mathcal{S}(\b x))\leq \cdots \leq \mathrm{dist}(\b \xi_N, \mathcal{S}(\b x))$. Then denote the ($k+1$)-th smallest distance value as $\mathrm{dist}^*(\b x) \coloneqq \mathrm{dist}(\b \xi_{k+1}, \mathcal{S}(\b x)) = s^*$. Combined with the analytical expression of the distance function in~\eqref{eq:dist_anal}, the pair $(\b x, \b r^*, s^*)$ satisfies the constraint \eqref{eq:MILP_reduced_part}. It remains to show the satisfaction of~\eqref{eq:MILP_reduced_q} given the pair $(\b x, \b r^*, s^*)$. Now consider two cases:\\
    1) If $\mathrm{dist}^*(\b x) = 0$, then $s^* = 0$. As~\eqref{eq:epsilonNs} must be satisfied for the given ($\b r^*$, $s^*$), this leads to a contradiction because $s^* = 0$ implies the non-positivity of the left-hand side of~\eqref{eq:epsilonNs}, but the right-hand side of~\eqref{eq:epsilonNs} is strictly positive due to the assumption $\theta > 0$. Therefore, we must have $\mathrm{dist}^*(\b x) > 0$.\\
    2) If $\mathrm{dist}^*(\b x) > 0$, then for any $p\in [P]$, when $i\geq k+1$, we have: $0 < s^*  = \mathrm{dist}^*(\b x) \leq \mathrm{dist}(\b \xi_i, \mathcal{S}(\b x))$. Given $0<\mathrm{dist}(\b \xi_i, \mathcal{S}(\b x))$, based on~\eqref{eq:dist_anal}, we have $\mathrm{dist}(\b \xi_i, \mathcal{S}(\b x)) = \min_{p\in [P]} \dfrac{\b b_p^\top \b \xi_i + d_p - \b a_p^\top \b x}{\| \b b_p \|_*}$, which leads to:
    $
       s^*  \leq  \min_{p\in [P]} \dfrac{\b b_p^\top \b \xi_i + d_p - \b a_p^\top \b x}{\| \b b_p \|_*}
       \leq  \dfrac{\b b_p^\top \b \xi_i + d_p - \b a_p^\top \b x}{\| \b b_p \|_*}.
    $
    In other words, for any $p \in [P]$, there will be at least $N-k$ elements in the set $\left\{\frac{\b b_p^\top \b \xi_i + d_p - \b a_p^\top \b x}{\| \b b_p \|_*} \right\}_{i\in [N]}$ not smaller than $s^*$. Since $q_p$ is defined as the $(k+1)$-th smallest value of $\{\b b_p^\top \b \xi_i\}_{i\in [N]}$, $\frac{q_p + d_p - \b a_p^\top \b x}{\| \b b_p \|_*}$ is the $(k+1)$-th smallest element in the set $\left\{\frac{\b b_p^\top \b \xi_i + d_p - \b a_p^\top \b x}{\| \b b_p \|_*} \right\}_{i\in [N]}$. Hence, it follows $s^* \leq \dfrac{q_p + d_p - \b a_p^\top \b x}{\| \b b_p \|_*}$, which leads to the satisfaction of~\eqref{eq:MILP_reduced_q}.
$\hfill\blacksquare$\endproof

\subsection*{Proof of Theorem~\ref{theo:propLP_sub_MILP}}

\proof{proof}
    By comparing $\mathcal{X}_\text{SFLA}(\b \kappa)$ in~\eqref{eq:SFLA} and $\mathcal{X}_\text{Exact}$ through its equivalent and strengthened reformulation in~\eqref{eq:MILP_constr_reduced}, the proof is apparent due to $\left(\dfrac{\b b_p^\top \b \xi_i + d_p - \b a_p^\top \b x}{\| \b b_p \|_*} \right)^+ \geq \kappa_i \left(\dfrac{\b b_p^\top \b \xi_i + d_p - \b a_p^\top \b x}{\| \b b_p \|_*} \right),\ \forall i \in [N]$ given $\kappa_i \in [0, 1]$.
$\hfill\blacksquare$\endproof





\subsection*{Proof of Theorem~\ref{theo:propLP_sup_LPori_subMILP}}

\proof{proof}
    It is equivalent to prove that for any $\b x \in \mathcal{X}_\text{LA}(\b\kappa)$, we have $\b x \in \mathcal{X}_\text{SFLA}(\b\kappa)$, namely that constraints~\eqref{eq:prop_LP_s_r}--\eqref{eq:prop_LP_q} can be satisfied.
    Based on Lemma \ref{lemma:our_exist}, for any $\b x \in \mathcal{X}_\text{LA}(\b\kappa)$, there always exists ($\b r^*$, $s^*$) such that $s^*$ is equal to the (k+1)-th smallest value amongst the set $\{ \widehat{\mathrm{dist}}(\b \xi_i, \mathcal{S}(\b x)) \}_{i\in [N]}$ and constraints~\eqref{eq:LP1_s_r}--\eqref{eq:LP1_last} are satisfied, which indicate the satisfaction of~\eqref{eq:prop_LP_s_r}--\eqref{eq:prop_LP_part}. 
    
    It remains to show the satisfaction of~\eqref{eq:prop_LP_q} given the same ($\b r^*$, $s^*$). Without loss of generality, suppose that we have the following ordering $\widehat{\mathrm{dist}}(\b \xi_1, \mathcal{S}(\b x)) \leq \cdots \leq \widehat{\mathrm{dist}}(\b \xi_N, \mathcal{S}(\b x))$. Then denote its (k+1)-th smallest element as $\widehat{\mathrm{dist}}^*(\b x) \coloneqq \widehat{\mathrm{dist}}(\b \xi_{k+1}, \mathcal{S}(\b x)) = s^*$. Based on~\eqref{eq:dist_approx}, we have $\widehat{\mathrm{dist}}^*(\b x) = \kappa_{k+1} \left(\min_{p\in [P]} \dfrac{\b b_p^\top \b \xi_{k+1} + d_p - \b a_p^\top \b x}{\| \b b_p \|_*}\right)$.
    Then for any $p\in [P]$, when $i\geq k+1$, we have $s^*  = \widehat{\mathrm{dist}}^*(\b x) \leq \widehat{\mathrm{dist}}(\b \xi_i, \mathcal{S}(\b x)) =  \kappa_i \left(\min_{p\in [P]} \dfrac{\b b_p^\top \b \xi_i + d_p - \b a_p^\top \b x}{\| \b b_p \|_*}\right)
           \leq  \kappa_i \left(\dfrac{\b b_p^\top \b \xi_i + d_p - \b a_p^\top \b x}{\| \b b_p \|_*} \right)$.
    In other words, for any $p \in [P]$, there will be at least $N-k$ elements in the set $\left\{\kappa_i (\b b_p^\top \b \xi_i + d_p - \b a_p^\top \b x) / \| \b b_p \|_* \right\}_{i\in [N]}$ not less than $s^*$. Since $q_p$ is defined as the $(k+1)$-th smallest value of $\left\{\b b_p^\top \b \xi_i\right\}_{i\in [N]}$, $(q_p + d_p - \b a_p^\top \b x) / \| \b b_p \|_*$ is the $(k+1)$-th smallest element of the set $\left\{ (\b b_p^\top \b \xi_i + d_p - \b a_p^\top \b x) / \| \b b_p \|_* \right\}_{i\in [N]}$. We can then deduce that $(q_p + d_p - \b a_p^\top \b x) / \| \b b_p \|_*$ is not less than the $(k+1)$-th smallest element of the set $\left\{\kappa_i (\b b_p^\top \b \xi_i + d_p - \b a_p^\top \b x) / \| \b b_p \|_* \right\}_{i\in [N]}$ due to $\kappa_i \in [0,1]$.
    Therefore, we must have $s^* \leq \dfrac{q_p + d_p - \b a_p^\top \b x}{\| \b b_p \|_*}$,
    which indicates the satisfaction of~\eqref{eq:prop_LP_q}.
$\hfill\blacksquare$\endproof

\subsection*{Proof of Corollary~\ref{cor:eq1}}

\proof{proof}
    Based on the definition of $[N]_p$, we have $\b b_p^\top \b \xi_i \geq q_p, \quad \forall i \in [N] \setminus [N]_p,\ p\in [P]$.
    Therefore, the constraint~\eqref{eq:prop_LP_q} in $\mathcal{X}_\text{SFLA}(\b 1)$ implies that for any $i \in [N] \setminus [N]_p$ with $p\in [P]$, we have $\dfrac{\b b_p^\top \b \xi_i + d_p - \b a_p^\top \b x}{\| \b b_p \|_*} \geq  \dfrac{q_p + d_p - \b a_p^\top \b x}{\| \b b_p \|_*} \geq s$.
    As $r_i \geq 0$, we further have $\frac{\b b_p^\top \b \xi_i + d_p - \b a_p^\top \b x}{\| \b b_p \|_*} \geq s \geq s - r_i$,
    which, combined with~\eqref{eq:prop_LP_part} in $\mathcal{X}_\text{SFLA}(\b 1)$, implies the satisfaction of~\eqref{eq:LP1_last} when $\b \kappa = \b 1$.
    The above process implies $\mathcal{X}_\text{SFLA}(\b 1) \subseteq \mathcal{X}_\text{LA}(\b 1)$. Combining with Theorem~\ref{theo:propLP_sup_LPori_subMILP}, we have $\mathcal{X}_\text{SFLA}(\b 1) = \mathcal{X}_\text{LA}(\b 1)$.
$\hfill\blacksquare$\endproof

\subsection*{Proof of Corollary \ref{cor:wcvar_la}}
\proof{proof}
In the W-CVaR formulation~\eqref{eq:CVaR_RHS}, constraint~\eqref{eq:wcvar-last} specifies the lower bound of $\beta$, $\beta \ge w_p \, \| \mathbf{b}_p \|_*,\forall p \in [P]$. Note that at optimality we can always safely set $\beta = \max_{p \in [P]} w_p \, \| \mathbf{b}_p \|_*$, which allows us to eliminate $\beta$ from the W-CVaR formulation~\eqref{eq:CVaR_RHS}. Next, applying the variable substitutions $r \leftarrow \frac{\b \alpha}{\beta}$, $-s \leftarrow \frac{\tau}{\beta}$, $\beta = \max_{p\in[P]}{w_p \| \b b_p \|_*}$ and after reorganizing, we obtain
\begin{subequations}\label{eq:CVaR_RHS_2}
    \begin{align}
        & s\geq 0, \b r \geq \b 0, \label{eq:CVaR_RHS_2a} \\
        & \epsilon N s - \sum\limits_{i\in [N]}r_i \geq \theta N,  \label{eq:CVaR_RHS_2b}\\
        & \frac{w_p \| \b b_p \|_*}{\max_{p\in[P]}{w_p \| \b b_p \|_*}}\left( \frac{\b b_p^\top \b \xi_i + d_p - \b a_p^\top \b x}{\| \b b_p \|_*} \right) \geq s-r_i, &&\hspace{-50pt} \forall i\in [N], p\in [P].  \label{eq:CVaR_RHS_2c}
    \end{align}
\end{subequations}
Comparing \eqref{eq:CVaR_RHS_2} and the LA formulation~\eqref{eq:LP_constr_1}, we see that the only difference is in constraint \eqref{eq:CVaR_RHS_2c} and constraint \eqref{eq:LP1_last}. Replacing $w_p$ with $w_p = \frac{1}{P}$, constraint~\eqref{eq:CVaR_RHS_2c} can be equivalently written as $\left( \dfrac{\b b_p^\top \b \xi_i + d_p - \b a_p^\top \b x}{\max_{p\in[P]}\| \b b_p \|_*} \right) \geq s-r_i,\ \forall i\in [N], p\in [P]$. This is a more conservative constraint compared to constraint \eqref{eq:LP1_last} in LA with $\b \kappa = \b 1$, because of 
$\left( \dfrac{\b b_p^\top \b \xi_i + d_p - \b a_p^\top \b x}{\| \b b_p \|_*} \right)\geq \left( \dfrac{\b b_p^\top \b \xi_i + d_p - \b a_p^\top \b x}{\max_{p\in[P]}\| \b b_p \|_*} \right) \geq s-r_i$. Therefore, we have $\mathcal{X}_\text{WCVaR}(\b 1 / P) \subseteq \mathcal{X}_\text{LA}(\b 1)$, and subsequently $\mathcal{X}_\text{WCVaR}(\b 1 / P)\subseteq \mathcal{X}_\text{SFLA}(\b 1)$ due to Corollary \ref{cor:eq1}. 
$\hfill\blacksquare$\endproof

\subsection*{Proof of Corollary~\ref{cor:obj_dependent_exact}}

\proof{proof}
    Proposition 9 in~\cite{chen2023approximations} demonstrates that the LA method can achieve an optimal value equivalent to the exact reformulation of RHS-WDRJCC by optimizing $\b\kappa$. We can then adjust the objective function $c(\cdot)$ and the set $\mathcal{X}$ to identify each individual feasible solution $\b x\in \mathcal{X}_\text{Exact}$, thereby establishing set equivalence. This set equivalence also applies to our proposed SFLA due to Theorem~\ref{theo:propLP_sup_LPori_subMILP}.
$\hfill\blacksquare$\endproof

\subsection*{Proof of Corollaries \ref{cor:wcvar} and \ref{cor:exact_conditions}}

\proof{proof}
The properties in Corollaries \ref{cor:wcvar} and \ref{cor:exact_conditions} hold for $\mathcal{X}_\text{LA}(\b \kappa)$ as was proved from the perspective of the optimal value in \cite{chen2023approximations}. This perspective can be extended to the whole set following our proof of Corollary~\ref{cor:obj_dependent_exact}, which in turn holds for our proposed SFLA $\mathcal{X}_\text{SFLA}(\b \kappa)$ due to Corollary \ref{cor:eq1}. In fact, the first condition in Corollary~\ref{cor:exact_conditions} can be relaxed from $\forall \b x \in \mathcal{X}_\text{Exact}$ to just the optimal solution $\b x^*$, as the primary interest is typically in the optimal solution rather than the entire set.
$\hfill\blacksquare$\endproof

\subsection*{Proof of Corollary \ref{cor:hyper}}

\proof{proof}
The proof is similar to that for Proposition 10 in \cite{chen2023approximations}.
$\hfill\blacksquare$\endproof

\subsection*{Proof of Theorem \ref{theo:equivalence_risk_radius}}

\proof{proof}
The exact reformulation of the WDRJCC is given by~\eqref{eq:wdrjcc_abstract}. We only need to focus on the constraint $\epsilon N s - \sum\limits_{i\in [N]}r_i \geq \theta N$, where we can show that for any feasible solution $(\b x, s, \b r)$ of one set, the same solution is also feasible to this constraint of the other set by having different $\theta$ or $\epsilon$ values. \\
$\mathcal{X}_\text{risk}(\theta_{\min}) \subseteq \mathcal{X}_\text{rad}(\epsilon_{\max})$: Let $(\b x^\text{risk}, s^\text{risk}, \b r^\text{risk})$ be any feasible solution to $\mathcal{X}_\text{risk}(\theta_{\min})$. We need to show the existence of $\theta$ such that $\epsilon_{\max} N s^\text{risk} - \sum\limits_{i\in [N]}r^\text{risk}_i \geq \theta N$ for some $\theta \geq \theta_{\min}$. Since $(\b x^\text{risk}, s^\text{risk}, \b r^\text{risk})$ is feasible to $\mathcal{X}_\text{risk}(\theta_{\min})$, we have $\epsilon N s^\text{risk} - \sum\limits_{i\in [N]}r^\text{risk}_i \geq \theta_{\min} N$ with $0\leq\epsilon \leq \epsilon_{\max}$. Then, it can be verified that setting $\theta = \theta_{\min}$ leads to the result due to $\epsilon_{\max} N s^\text{risk} - \sum\limits_{i\in [N]}r^\text{risk}_i \geq \epsilon N s^\text{risk} - \sum\limits_{i\in [N]}r^\text{risk}_i \geq \theta_{\min} N $.
\\[1.1em]
$\mathcal{X}_\text{rad}(\epsilon_{\max}) \subseteq\mathcal{X}_\text{risk}(\theta_{\min})$: Denote by $(\b x^\text{rad}, s^\text{rad}, \b r^\text{rad})$ be any feasible solution to $\mathcal{X}_\text{rad}(\theta_{\min})$. We need to show the existence of $\epsilon$ such that $ \epsilon N s^\text{rad} - \sum\limits_{i\in [N]}r^\text{rad}_i \geq \theta_{\min} N$ for some $ \epsilon \in [0, \epsilon_{\max}]$. Again, because $(\b x^\text{rad}, s^\text{rad}, \b r^\text{rad})$ is feasible to $\mathcal{X}_\text{rad}(\epsilon_{\max})$, we have $\epsilon_{\max} N s^\text{rad} - \sum\limits_{i\in [N]}r^\text{rad}_i \geq \theta N$ with $\theta \geq \theta_{\min}$. We can then set $\epsilon = \epsilon_{\max}$ to obtain the desired inequality, because of $\epsilon_{\max} N s^\text{rad} - \sum\limits_{i\in [N]}r^\text{rad}_i \geq \theta N \geq \theta_{\min} N$.
$\hfill\blacksquare$\endproof

}

\section{Supplemental Information for the Unit Commitment Problems}

\begin{figure}
    \centering 
    \vspace{-3.5 mm}
    \includegraphics[width=.49\linewidth]{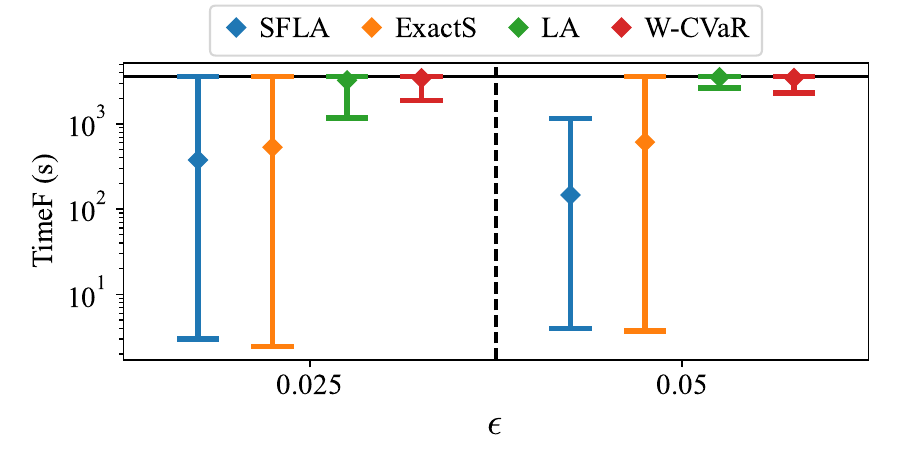}
    \includegraphics[width=.49\linewidth]{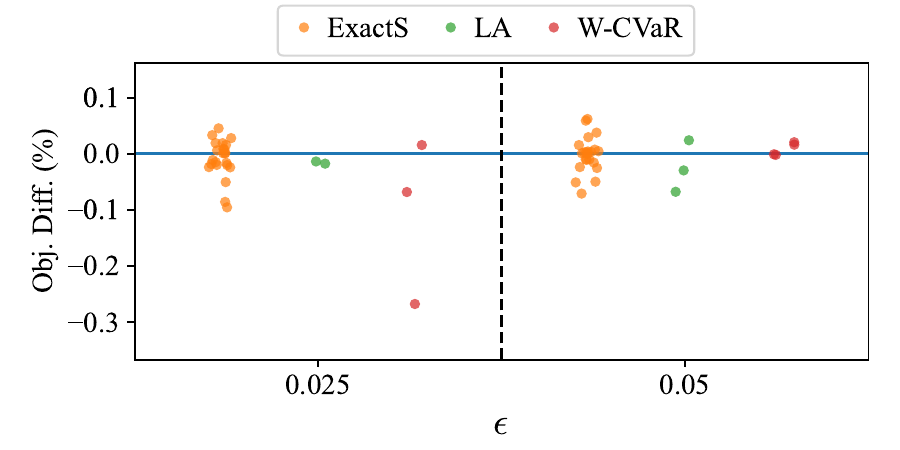}
    \begin{minipage}{\textwidth}
    \hspace{1em}
    \begin{minipage}{\textwidth}
        \centering
        \vspace{-3em}  
        \subfloat[\label{fig:suc_allmethods_timeF}]{\phantom{\rule{0.49\textwidth}{0pt}}}
        \subfloat[\label{fig:suc_allmethods_gap}]{\phantom{\rule{0.49\textwidth}{0pt}}}
    \end{minipage}
    \end{minipage}
    \vspace{-3.3em}
    \caption{ 
    Performance comparison of the proposed SFLA with benchmarks for the UC problem: ExactS, LA, W-CVaR, and Bonferroni for the UC problem.}
    \parbox{\textwidth}{\vspace{2mm}
    \footnotesize{\textit{Note.} The Wasserstein radius $\theta$ is set to $0.5$ and the optimization horizon is set to $T=24$. (a) Comparison of \emph{TimeF} (s), where dots represent the mean value of the 30 random runs, with error bars indicating the 95\% percentile interval. The black horizontal line represents the one-hour \texttt{TimeLimit}. (b) Comparison of \emph{Obj. Diff.} (\%), where each dot represents the result for one of the 30 runs. The blue horizontal line indicates zero difference between the optimal value of a benchmark method and the proposed SFLA. Values higher than this horizontal line indicate that the benchmark achieves a lower cost (better optimality) than the proposed SFLA. Note that we only keep certain random runs where the method to be compared and the proposed method are both solved to optimum. Many runs of LA and W-CVaR fail to converge to the optimum, and thus the plot only shows a few or zero dots for these two methods. The Bonferroni approximation is infeasible due to its over-conservativeness for all runs and is therefore not displayed. }
    \vspace{3mm}}
    \label{fig:suc_allmethods}
\end{figure}

\begin{figure}
     \FIGURE
    {\includegraphics[width=1\linewidth]{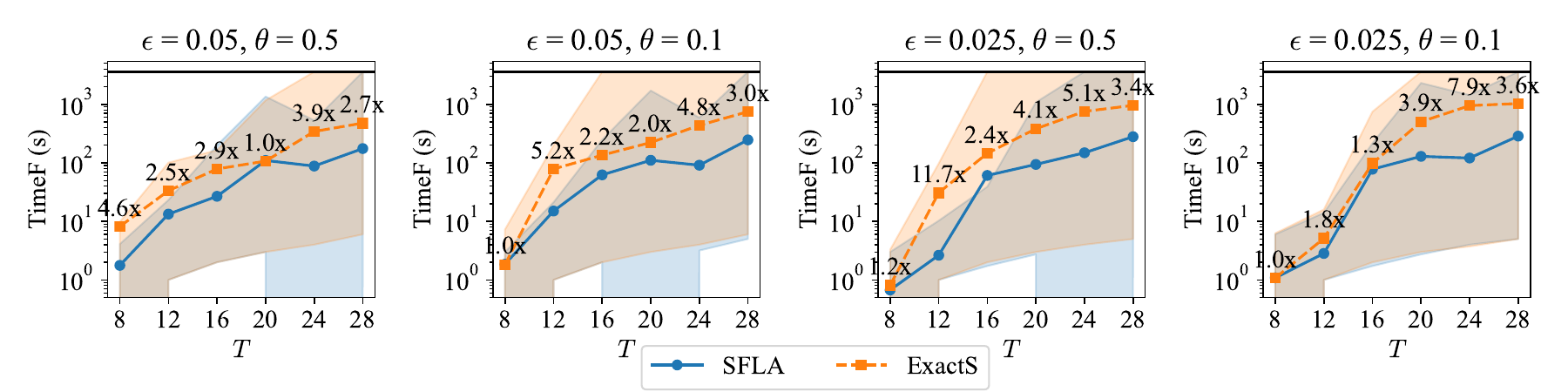}\vspace{-2mm}}
    {Computation time to obtain the first comparable high-quality solution under different optimization horizons $T$ for the proposed SFLA and the benchmark ExactS for the UC problem. \label{fig:suc_timeF_allT}}
    {Dots represent the mean value of the 150 random runs, with shaded areas indicating the 95\% percentile interval (from 2.5th to 97.5th). The black horizontal line represents the one-hour \texttt{TimeLimit}. The number of historical samples is set to $N=100$. Numbers ending with ``x'' represent the speedup in the average computing time of the proposed SFLA compared to ExactS.}
\end{figure}
\begin{figure}
     \FIGURE
    {\includegraphics[width=1\linewidth]{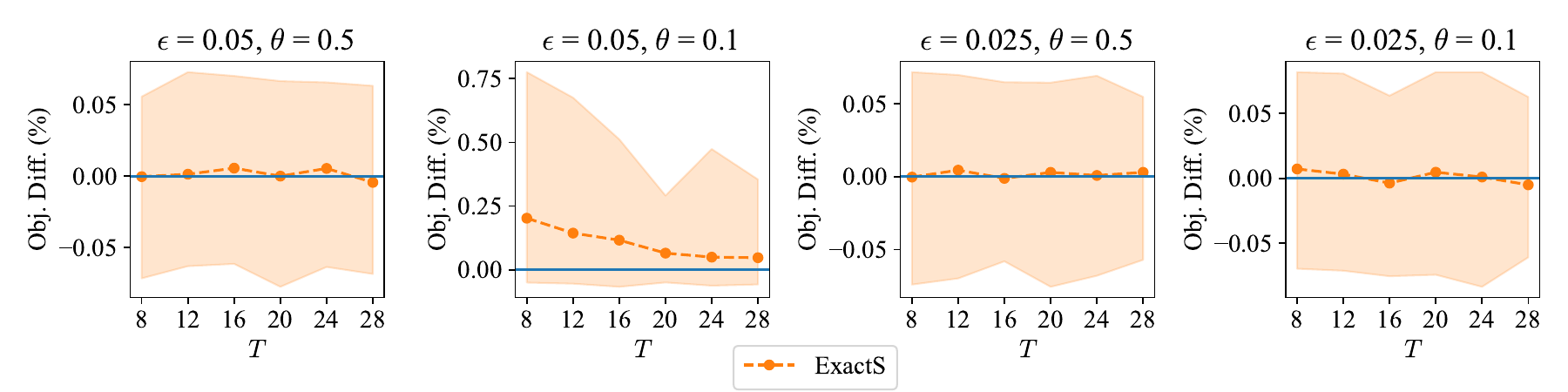}\vspace{-2mm}}
    {Comparison of optimality for different optimization horizons $T$ for the proposed SFLA and the benchmark ExactS for the UC problem. \label{fig:suc_gap_allT}}
    {The round dots represent the mean value of the 150 random runs, with error bars indicating the 95\% percentile interval. The blue horizontal line indicates a zero difference of the optimal value of a benchmark method compared to the proposed SFLA. Values higher than this horizontal line indicate that the benchmark achieves lower cost (better optimality) than the proposed SFLA. The number of historical samples is set to $N=100$.}
\end{figure}

\subsection{Case Study Settings for the Unit Commitment Problem}
\label{appen:setting_UC}

To ensure robust conclusions, we conduct $30$ random runs for each parameter setting of $\epsilon$, $\theta$, $N$, and $T$, increasing the number of random runs to $150$ for specific evaluations when needed. In each random run, we:

\noindent 1) Sample a distinct set of generator cost parameters from a uniform distribution based on the range suggested by~\cite{xu2017application}, which impacts the objective function of the UC problem.

\noindent 2) Sample a different set of historical samples $\{ \b \xi_i \}_{i \in [N]}$ for the random wind forecasting error, uniformly from a total of $1000$ training data samples, which affects the feasible region of the UC problem.

\noindent 3) Select a random starting simulation hour, uniformly sampled from the $24$ hours of a day, which affects the demand and generation profiles and thereby alters the feasible region of the UC problem.

\noindent 4) Set a different seed for the optimization solver, affecting the solver behaviors such as tie-breaking rules.

\begin{figure}
     \FIGURE
    {\hfill
    \subfloat[ ]{\includegraphics[width=.5\linewidth]{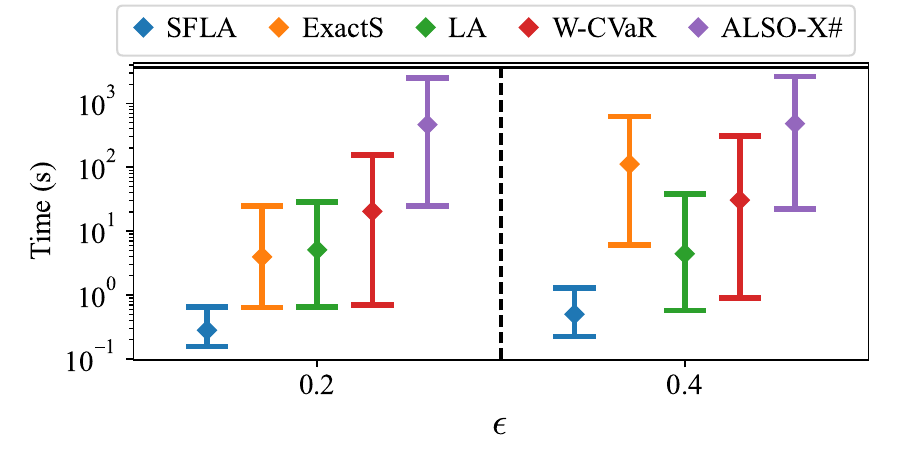}\label{fig:suc_also_allmethods_timeF}}\hfill
    \subfloat[ ]{\includegraphics[width=.5\linewidth]{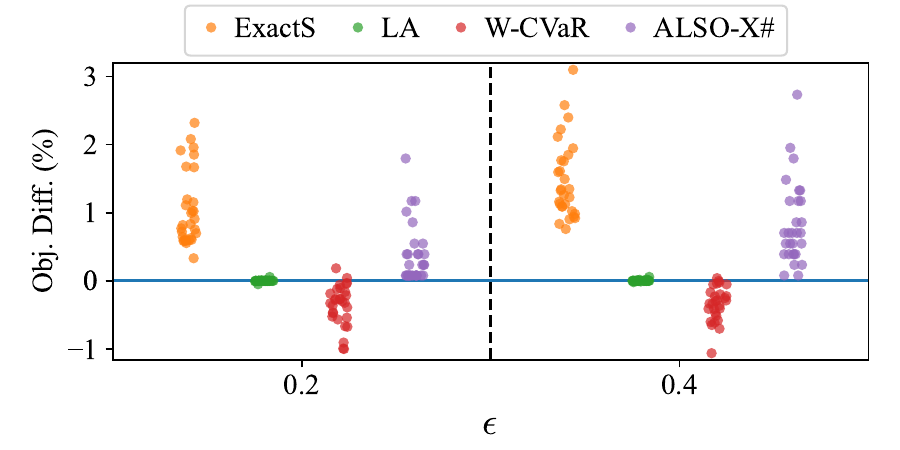}\label{fig:suc_also_allmethods_gap}}
    \hfill} 
{Performance comparison of the proposed SFLA with benchmarks for the UC problem: ExactS, LA, W-CVaR, and ALSO-X\# for the UC problem. \label{fig:suc_also_allmethods}}
{The Wasserstein radius $\theta$ is set to $0.1$ and the optimization horizon is set to $T=2$ as otherwise ALSO-X\# is unsolvable. The number of generators also reduces to $10$ for tractability. All graphical elements are defined as in Figure~\ref{fig:suc_allmethods} and are not repeated here.
}
\end{figure}

\subsection{Supplemental Numerical Experiments of the Unit Commitment Problem}
\label{appen:sup_res_UC}

Figure~\ref{fig:suc_allmethods_timeF} compares all methods with an optimization horizon $T=24$. As can be seen, in terms of mean values, the proposed SFLA can have $4 \times$ faster UC computation than the exact reformulation (ExactS) and even $20 \times$ faster than the two linear approximation schemes, namely LA and W-CVaR. This demonstrates the computational efficiency of the proposed SFLA, and also the effectiveness of the strengthening process that makes ExactS even more computationally efficient than linear approximations. Figure~\ref{fig:suc_allmethods_gap} shows that the proposed SFLA has optimality comparable to ExactS and LA, where the minor difference of less than $0.1\%$ is attributed to the preset $0.1\%$ \texttt{MIPGap}. We also observed a random run of W-CVaR with a 0.25\% worse optimality compared to the proposed SFLA, suggesting the higher conservativeness of W-CVaR. The Bonferroni approximation is infeasible in all cases due to its over-conservatism.
Overall, the observations above support both the computational gain and the approximation quality of the proposed SFLA.

Figures~\ref{fig:suc_timeF_allT} and~\ref{fig:suc_gap_allT} compare the proposed SFLA and ExactS under progressively increasing optimization horizons $T$, a primary driver of problem complexity. Figure~\ref{fig:suc_timeF_allT} shows that the proposed SFLA sustains its computational speedup as problem complexity increases. In contrast, the optimality differences illustrated in Figure~\ref{fig:suc_gap_allT} exhibit a flat or even decreasing trend for $(\epsilon=0.05, \theta=0.1)$. These results indicate that the proposed SFLA preserves computational efficiency in more complex cases without sacrificing optimality, making it particularly advantageous for practical high-complexity applications such as UC.

We also compared against the ALSO-X\# algorithm~\citep{jiang2024also}. This algorithm interprets W-CVaR as an iterative bisection process. In each iteration, a new bound of the objective (determined by bisection) is added as a constraint, and the process terminates once the W-CVaR evaluation (a conservative approximation of WDRJCC) is satisfied. ALSO-X\# modifies this procedure by replacing the conservative W-CVaR stopping criterion with the exact WDRJCC, leveraging the fact that checking WDRJCC feasibility is tractable. This substitution reduces conservativeness and can even recover exact optimality under certain conditions, since the stopping condition becomes less restrictive. However, in each iteration, an optimization problem of the same scale as W-CVaR must be solved, leading to computation times that are several times longer than those of W-CVaR, as reported in their numerical experiments~\citep{jiang2024also}. Moreover, ALSO-X\# requires initial upper and lower bounds of the objective as inputs, and the quality of these bounds can significantly affect the efficiency of the subsequent bisection process. In fact,~\cite{jiang2024also} used W-CVaR to calculate the upper bound. Since Figure~\ref{fig:suc_allmethods} demonstrates that W-CVaR is both slower and can yield worse optimality compared to SFLA, it follows that ALSO-X\# will inevitably be slower than SFLA in obtaining high-quality solutions when W-CVaR is used as the starting point.

Since ALSO-X\# only requires initial bounds of the objective as inputs, it can be viewed as an add-on that further improves the solution quality of the proposed SFLA. Accordingly, in our experiments, we use the SFLA output as the initial upper bound of ALSO-X\#, set the initial lower bound to 95\% of this upper bound, and record only the bisection time as the computation time of ALSO-X\#.
As illustrated in Figure~\ref{fig:suc_also_allmethods}, ALSO-X\# is still slower than the MIP-based ExactS method. The main reason is that ALSO-X\# does not exploit the strengthening opportunities available in RHS-WDRJCC (a limitation similar to W-CVaR). A promising direction for future work may be the combination of SFLA with ALSO-X\#. Figure~\ref{fig:suc_also_allmethods_gap} shows that while ALSO-X\# improves solution quality, it still underperforms compared to the exact method. Notably, in this case, the exact method achieves about $3$\% better optimality than the proposed SFLA because we set a relatively high risk level $\epsilon$ to emphasize the improvement brought by ALSO-X\#.

\section{Supplemental Information for the Bilevel Bidding Problem}
\label{appen:sup_res_bilevel}

\subsection{Case Study Settings for the Bilevel Bidding Problem}
\label{appen:setting_bilevel}

To ensure reliable conclusions, we conduct $30$ random runs for each parameter setting. In each random run, we:

\noindent 1) Sample a distinct set of generator bid prices from a uniform distribution ranging between $80\%$ and $120\%$ of the bid prices used by~\cite{8036231}, which influences the objective function of the bilevel problem.

\noindent 2) Sample a different set of historical samples $\{ \b \xi_i \}_{i \in [N]}$ for the random wind forecasting error, uniformly from a total of $1000$ training data samples, which affects the feasible region of the bilevel problem.

\noindent 3) Select a random starting simulation hour, uniformly sampled from the $24$ hours of a day, which affects the demand and generation profiles and thereby alters the feasible region of the bilevel problem.

\noindent 4) Set a different seed for the optimization solver, affecting the solver behaviors such as tie-breaking rules.

\subsection{Supplemental Evaluation Metrics}

The exclusive evaluation metrics used in this supplemental case study include:

\noindent 1) \emph{N Infeasible}: This is the number of random runs that are infeasible. We introduce this metric because the high complexity of the bilevel problem introduces numerical issues that may make solvers incorrectly assess a feasible instance as infeasible. 
A large value of this metric indicates that the method is prone to numerical issues, especially for LA and W-CVaR.

\noindent 2) \emph{Bilevel Profit} (k\$): This represents the optimal profit of storage systems, which is the objective obtained by solving the bilevel problem.
 
\noindent 3) \emph{Actual Profit} (k\$): The market-clearing process is a single-level problem, where the market operator can use the exact MIP reformulation of RHS-WDRJCC to maximize social welfare. Therefore, the use of ExactS may cause \emph{Bilevel Profit} to differ from \emph{Actual Profit}. We calculate \emph{Actual Profit} through the following steps: i) Retrieve the bids and offers from the bilevel problem solutions; ii) Adjust bid prices slightly downward and offer prices slightly upward (by \$\,$10^{-5}$), setting bid quantities to zero when the bilevel-cleared quantities are zero---these adjustments ensure the desired acceptance or rejection behaviors modeled in the bilevel problem; iii) Input these modified bids and offers into the actual market clearing problem with ExactS (MIP); iv) Calculate the actual market clearing price, which is then used to compute the \emph{Actual Profit} of the storage operator. Due to the MIP formulation in ExactS, optimal dual variables cannot be directly used to determine marginal prices. Instead, we use sensitivity analysis to calculate the clearing price, by comparing the change of the objective function of the actual market clearing after a $10^{-3}$ MW increase in demand for each time step.

\noindent 4) \emph{Actual Profit Diff} (k\$): The difference of \emph{Actual Profit} between a benchmark method and that calculated by the proposed SFLA. A negative value represents a lower profit compared to SFLA. We do not use the percentage difference as is in our UC case because the profit could be close to zero, leading to exacerbated results. When calculating \emph{Actual Profit Diff}, we consider only those random runs where both the proposed SFLA and the benchmark identify feasible solutions within the prescribed \texttt{TimeLimit}.

\begin{figure}
    \centering    
    \includegraphics[width=\columnwidth]{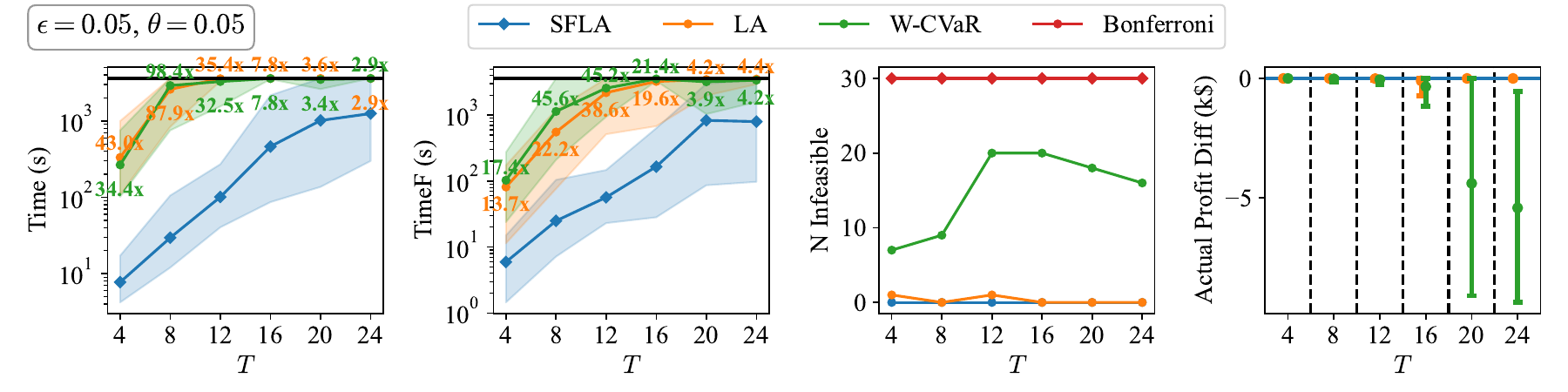}
    \begin{minipage}{\textwidth}
    \hspace{1em}
    \begin{minipage}{\textwidth}
        \centering
        \vspace{-3em}  
        \subfloat[]{\phantom{\rule{0.247\textwidth}{0pt}}}
        \subfloat[]{\phantom{\rule{0.247\textwidth}{0pt}}}
        \subfloat[\label{fig:n_infeasible_eps0.05_theta0.05}]{\phantom{\rule{0.247\textwidth}{0pt}}}
        \subfloat[\label{fig:profit_diff_eps0.05_theta0.05}]{\phantom{\rule{0.247\textwidth}{0pt}}}
    \end{minipage}
    \end{minipage}
    \begin{minipage}{\textwidth}
    \vspace{-2mm}
    \includegraphics[width=\columnwidth]{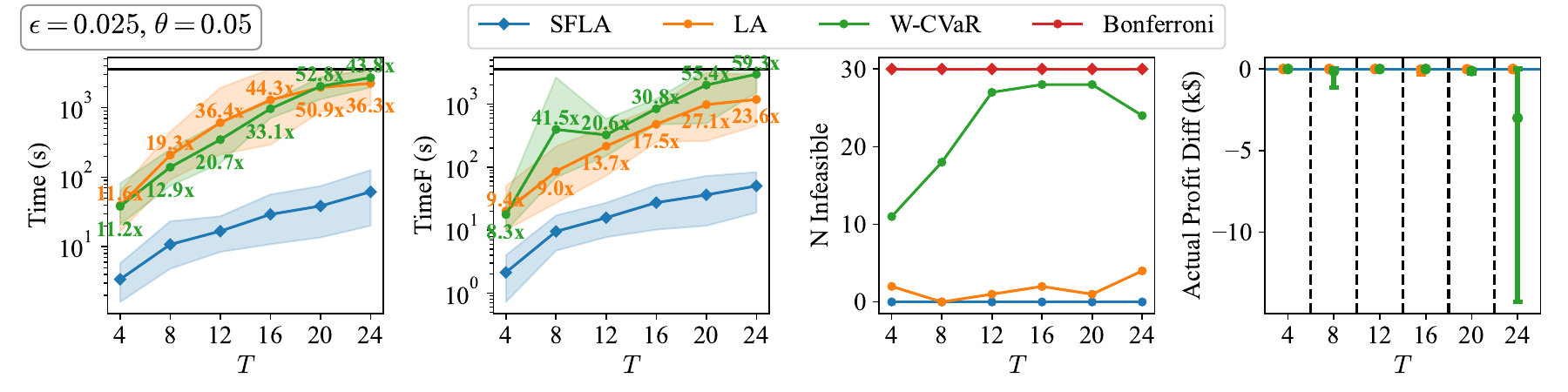}
    \end{minipage}
    \begin{minipage}{\textwidth}
    \hspace{1em}
    \begin{minipage}{\textwidth}
        \vspace{-1.5em}  
        \subfloat[]{\phantom{\rule{0.247\textwidth}{0pt}}}
        \subfloat[]{\phantom{\rule{0.247\textwidth}{0pt}}}
        \subfloat[\label{fig:n_infeasible_eps0.025_theta0.05}]{\phantom{\rule{0.247\textwidth}{0pt}}}
        \subfloat[\label{fig:profit_diff_eps0.025_theta0.05}]{\phantom{\rule{0.247\textwidth}{0pt}}}
    \end{minipage}
    \end{minipage}
    \vspace{-2.3em}
    \caption{ 
    Supplemental performance comparison for the bilevel problem.}
    \parbox{\textwidth}{\vspace{-4mm}
    \footnotesize{\textit{Note.} The results supplement Figure \ref{fig:bilevel_results} for $\theta=0.05$.}
    }
    \label{fig:bilevel_results_sup}
\end{figure}

\subsection{Supplemental Analysis for the Bilevel Bidding Problem}\label{appen:sup_bilevel}

This section provides supplemental analysis on \emph{N Infeasible} and the actual profits that can be attained by using the proposed SFLA and other benchmarks.
As can be observed in the third column of Figures~\ref{fig:bilevel_results} and \ref{fig:bilevel_results_sup}, W-CVaR shows a large number of infeasible runs. These infeasibilities are resolved by setting \texttt{NumericalFocus=3}, indicating numerical issues. While LA is more numerically stable, it still incurs $22$ infeasible runs in total. In contrast, the proposed SFLA has only one infeasible run, demonstrating its higher numerical stability. Note that the infeasibility of the Bonferroni approximation is caused by its over-conservativeness that cannot be recovered by setting \texttt{NumericalFocus=3}.

The fourth column of Figures~\ref{fig:bilevel_results} and \ref{fig:bilevel_results_sup} compares the \emph{Actual Profit} of the benchmarks with that of the proposed SFLA. Among the feasible cases, Bonferroni approximation leads to significantly lower profit than other methods because of its over-conservativeness. It is noteworthy that both LA and W-CVaR can occasionally yield either higher or lower profits than the proposed SFLA. This variability arises because the market-clearing model in the bilevel problem does not precisely replicate the actual market-clearing process (cleared with exact MIP-based RHS-WDRJCC), and bilevel problems admit multiple solutions even at optimality. Thus, even when different methods yield the same \emph{Bilevel Profit}, discrepancies can still arise in \emph{Actual Profit}. It should be noted that W-CVaR shows significantly lower \emph{Actual Profit}, especially for a large $T$, because W-CVaR fails to find a good feasible solution within the \texttt{TimeLimit}, rather than due to over-conservativeness.

Finally, among all $720$ runs displayed in Figures~\ref{fig:bilevel_results} and \ref{fig:bilevel_results_sup}, there are only three instances where \emph{Bilevel Profit} of the proposed SFLA is slightly (less than $0.009$ k\$) below that of LA and W-CVaR. This is caused by another aspect of numerical issue of the proposed SFLA that leads to an incorrect upper bound proved. Setting \texttt{NumericalFocus=3} of \texttt{Gurobi} effectively mitigates this issue. Compared to the numerical issue that renders the problem infeasible in LA and W-CVaR, the numerical issue associated with the proposed SFLA is minor.

\begin{figure}
     \FIGURE
    {\includegraphics[width=.8\linewidth]{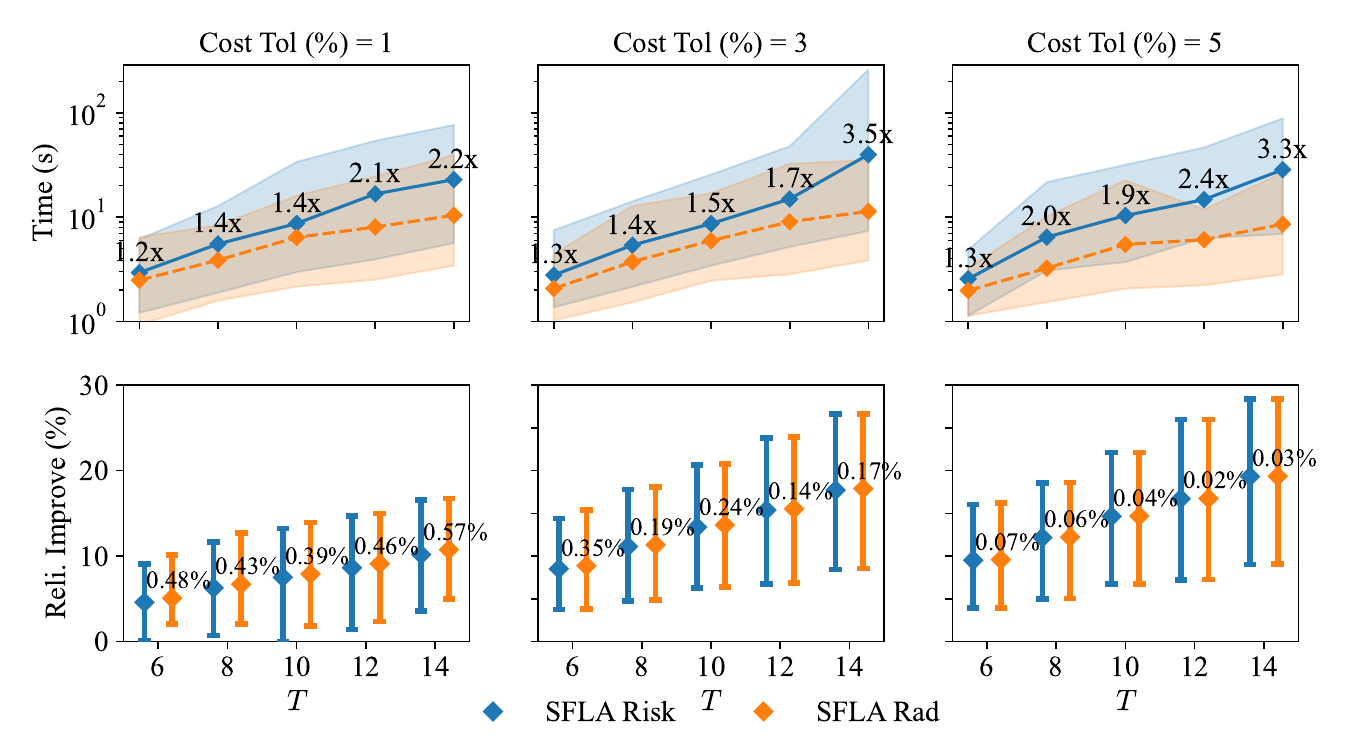} \vspace{-3mm}} 
    {Performance comparison of risk minimization and radius maximization using the proposed SFLA. \label{fig:suc_epstheta_evst}}
    {
    The square dots represent the mean values over 30 random runs, with shaded areas and error bars indicating the 95\% percentile interval. $T$ is the length of the optimization horizon. In the first row, the numbers ending with ``x'' are the speed-up of the average computing time of SFLA Rad compared to SFLA Risk. In the second row, the numbers represent the increased absolute amount of ``Reli. Improvement'' of SFLA Rad compared to SFLA Risk.
    }
\end{figure}



\section{Supplemental Results for Risk Minimization vs Radius Maximization} \label{res:risk_vs_radius}

Figure~\ref{fig:suc_epstheta_evst} compares the proposed SFLA in risk minimization (SFLA Risk) and radius maximization (SFLA Rad). The first row shows that radius maximization is more computationally tractable and can lead to up to 3.5x computational speedup compared to risk minimization. In the second row, the result is surprising as radius maximization can even be more robust than risk minimization, and this ``greater robustness'' result is more significant when the cost reduction budget is low (\emph{Cost Tol}=1\%). 


\end{APPENDIX}

\end{document}